\RequirePackage{ifpdf}
\ifpdf 
\documentclass[pdftex]{sigma}
\else
\documentclass{sigma}
\fi
\numberwithin{equation}{section}
\numberwithin{theorem}{section}
\numberwithin{lemma}{section}

\usepackage{enumerate}

\begin{document}

\allowdisplaybreaks

\renewcommand{\thefootnote}{$\star$}

\renewcommand{\PaperNumber}{066}

\FirstPageHeading

\ShortArticleName{Holonomy and Projective Equivalence in 4-Dimensional Lorentz Manifolds}

\ArticleName{Holonomy and Projective Equivalence \\ in 4-Dimensional Lorentz Manifolds\footnote{This paper is a
contribution to the Special Issue ``\'Elie Cartan and Dif\/ferential Geometry''. The
full collection is available at
\href{http://www.emis.de/journals/SIGMA/Cartan.html}{http://www.emis.de/journals/SIGMA/Cartan.html}}}

\Author{Graham S. HALL~$^\dag$ and David P. LONIE~$^\ddag$}

\AuthorNameForHeading{G.S. Hall and D.P. Lonie}

\Address{$^\dag$~Department of Mathematical Sciences, University of Aberdeen,\\
\hphantom{$^\dag$}~Meston Building, Aberdeen, AB24 3UE, Scotland, UK}
\EmailD{\href{mailto:g.hall@abdn.ac.uk}{g.hall@abdn.ac.uk}}

\Address{$^\ddag$~108e Anderson Drive, Aberdeen, AB15 6BW, Scotland, UK}
\EmailD{\href{mailto:DLonie@aol.com}{DLonie@aol.com}}

\ArticleDates{Received March 18, 2009, in f\/inal form June 11, 2009;  Published online June 29, 2009}

\Abstract{A study is made of 4-dimensional Lorentz manifolds which
are projectively related, that is, whose Levi-Civita connections
give rise to the same (unparameterised) geodesics. A brief review of
some relevant recent work is provided and a list of new results
connecting projective relatedness and the holonomy type of the
Lorentz manifold in question is given. This necessitates a review of
the possible holonomy groups for such manifolds which, in turn,
requires a certain convenient classif\/ication of the associated
curvature tensors. These reviews are provided.}

\Keywords{projective structure; holonomy; Lorentz manifolds; geodesic equivalence}

\Classification{53C29; 53C22; 53C50}

\rightline{\it This paper is dedicated to the memory of \'{E}lie Cartan}

\section{Introduction}

\'{E}lie Cartan (1869--1951) was one of the world's leading geometers
and it is to his memory that the authors dedicate this paper. One of
Cartan's main interests lay in the crucially important study of
connections on manifolds and its applications to theoretical
physics. The present paper will proceed in a similar vein by
presenting a discussion of holonomy theory on a 4-dimensional
manifold which admits a Lorentz metric. This leads naturally to
applications in Einstein's general theory of relativity. Here the
application will be to the projective structure of space-times and
which itself is closely related to the principle of equivalence in
Einstein's theory. A few weeks before this paper was begun, the
authors learned, sadly, of the death of \'{E}lie Cartan's son,
Henri, at the age of 104. Henri Cartan was also a world leader in
geometry and this father and son combination laid down the
foundations for a great deal of the research presently undertaken in
dif\/ferential geometry.

This paper will be arranged in the following way. Section~\ref{section2} will be
used to introduce general notation and in Section~\ref{section3} a classif\/ication
of the curvature tensor will be introduced and which will prove
useful in what is to follow. Also included in Section~\ref{section3} is a
discussion of certain relationships between the metric, connection
and curvature structures on a space-time. In Section~\ref{section4} a review of
holonomy theory will be given and, in particular, as it applies to
4-dimensional Lorentz manifolds. In Section~\ref{section5}, a discussion of
projective structure will be presented. This will be followed, in
Sections~\ref{section6} and~\ref{section7}, with several theorems which show a tight
relationship between projective relatedness and holonomy type. A
brief summary of the paper is given in Section~\ref{section8}.

\section{Notation and preliminary remarks}\label{section2}

Throughout this paper, $M$ will denote a 4-dimensional, Hausdorf\/f,
connected, smooth manifold which admits a smooth metric $g$ of
Lorentz signature $(-,+,+,+)$. The pair $(M,g)$ is called a~{\em
space-time}. It follows that the (usual manifold) topology on $M$ is
necessarily second coun\-tab\-le~\cite{1} and hence paracompact. All
structures on $M$  will be assumed smooth (where this is sensible).
The unique symmetric Levi-Civita connection arising on $M$ through
$g$ is denoted by~$\nabla$ and, in a coordinate domain of~$M$, its
Christof\/fel symbols are written $\Gamma^a_{bc}$. The type $(1,3)$
curvature tensor associated with $\nabla$ is denoted by ${\rm Riem}$ and
its (coordinate) components are written $R^a_{\ bcd}$. The Ricci
tensor, ${\rm Ricc}$, derived from ${\rm Riem}$, has components $R_{ab}=R^c_{\
acb}$ and $R=R_{ab}g^{ab}$ is the Ricci scalar. The Weyl type
$(1,3)$ {\em conformal tensor} $C$ has components $C^a_{\ bcd}$
given by
\begin{equation}\label{E1}
C^a_{\ bcd}=R^a_{\ bcd}-E^a_{\ bcd}-\tfrac{R}{12}\big(\delta^a_{\ c}g_{bd}-\delta^a_{\ d}g_{bc}\big),
\end{equation}
where $E$ is the tensor with components
\begin{equation}\label{E2}
E^a_{\ bcd}=\tfrac{1}{12}\big(\tilde{R}^a_{\ c}g_{bd}-\tilde{R}^a_{\ d}g_{bc}+\delta^a_{\ c}\tilde{R}_{bd}-\delta^a_{\ d}\tilde{R}_{bc}\big),
\end{equation}
and where $\widetilde{\rm Ricc}$ is the trace-free Ricci tensor with
components $\widetilde{R}_{ab}=R_{ab}-\tfrac{R}{4}g_{ab}$.  At any $m\in
M$, the tensors $E$ and $\widetilde{\rm Ricc}$ uniquely (algebraically)
determine each other and
$E=0\Leftrightarrow\widetilde{\rm Ricc}=0\Leftrightarrow$ the {\em Einstein
space} condition holding at $m$.

For $m\in M$, $T_mM$ denotes the {\em tangent space to $M$ at $m$}
and this will,  for convenience, be identif\/ied with the cotangent
space $T^*_mM$ to $M$ at $m$, through the metric $g(m)$ by index
raising and lowering. Thus the liberty will be taken of using the
same symbol for members of $T_mM$ and $T^*_mM$ which are so related
and similarly for other tensor spaces. A tetrad (that is, a basis
for $T_{m}M$) $u,x,y,z\in T_mM$ is called {\em orthonormal} if and
only if the only non-vanishing inner products between tetrad members
are $-g(u,u)=g(x,x)=g(y,y)=g(z,z)=1$ and a tetrad $l,n,x,y\in T_mM$
is called {\em null} if and only if the only non-vanishing inner
products between tetrad members are $g(l,n)=g(x,x)=g(y,y)=1$. In
this case, $l$ and $n$ are null vectors. Another condition, the {\em
non-f\/lat condition}, will be imposed on $(M,g)$, meaning that $\rm Riem$
does not vanish over any non-empty open subset of $M$. This is a
physical requirement and is there to prevent gravitational shielding
in general relativity theory.

Let $\Lambda_mM$ denote the 6-dimensional vector space of all tensor
type $(2,0)$ 2-forms at $m$.  This vector space can be associated,
using the metric $g(m)$, with the vector spaces of tensor type
$(0,2)$ 2-forms at $m$ and of type $(1,1)$ tensors at $m$ which are
skew-self adjoint with respect to~$g(m)$, through the component
identif\/ications
\[ (F\in\Lambda_mM$) $F^{ab}(=-F^{ba})\rightarrow
F_{ab}\equiv g_{ac}g_{bd}F^{cd}\rightarrow F^a_{\ b}=g_{cb}F^{ac}.\]
Any member $F$ of any of these vector spaces will be referred to as
a {\em bivector} (at $m$) and written symbolically as
$F\in\Lambda_mM$. Any $F\in\Lambda_mM$, $F\neq0$, has (matrix) rank
2 or 4. If~$F$ has rank~2 it is called {\em simple} and may be
written in components as $F^{ab}=p^aq^b-q^ap^b$ for $p,q\in T_mM$.
Although~$F$ does not determine $p$ and $q$ it does uniquely
determine the 2-dimensional subspace (referred to as a {\em
$2$-space}) of $T_mM$ spanned by $p$ and $q$ and which is called the
{\em blade} of $F$. A simple bivector~$F$ at $m$ is then called {\em
timelike}, {\em spacelike} or {\em null} according as its blade is,
respectively, a timelike, spacelike or null 2-space at $m$. If
$F\in\Lambda_mM$ has rank 4 it is called {\em non-simple} and may be
written as $F=G+H$ where $G$ and $H$ are simple bivectors with $G$
timelike and $H$ spacelike and where the blades of $G$ and $H$ are
uniquely determined by $F$ and are orthogonal complements of each
other. They will be collectively called the {\em canonical pair of
blades} of $F$. The Hodge duality operator on bivectors is denoted
by $*$. Then $F$ is simple if and only if $\overset{*}{F}$ is simple
and the blades of $F$ and $\overset{*}{F}$ are orthogonal
complements of each other. In this case $F$ is spacelike,
respectively, timelike or null, if and only if $\overset{*}{F}$ is
timelike, respectively, spacelike or null. If $F$ is simple and
either spacelike or timelike, the blades of $F$ and $\overset{*}{F}$
are complementary (that is, their union spans $T_mM$). This is not
the case if $F$ is null. In the above expression for a non-simple
bivector $F$, $G$ and $H$ are (multiples of) duals of each other.
For any bivector $F$ at $m$, $F$ and $\overset{*}{F}$ are
independent bivectors at $m$ and $\overset{**}{F}=-F$. If $F$ is
a~simple bivector at $m$ with $F^{ab}=p^aq^b-q^ap^b$, $p,q\in T_mM$,
then $F$ (or its blade) is sometimes written as $p\wedge q$. It is
convenient, on occasions, to use round and square brackets to denote
the usual symmetrisation and skew-symmetrisation of indices,
respectively. As a general remark on notation, both coordinate and
coordinate-free notation will be used, depending on relative
convenience. More details on these aspects of Minkowski geometry may
be found, for example, in~\cite{2}.

\section{Curvature structure of space-times}\label{section3}

Because of the algebraic symmetries of ${\rm Riem}$, one may introduce the
{\em curvature maps}  $f$ and $\tilde{f}$ from the vector space of
bivectors to itself at $m$ (recalling the liberties taken with this
vector space in Section~\ref{section2}) by
\begin{equation}\label{E3}
f: \ F^{ab}\rightarrow R^{a}_{\ bcd}F^{cd}, \qquad \tilde{f}: \
F^{ab}\rightarrow R^{ab}_{\ \ cd}F^{cd}.
\end{equation}
The maps $f$ and $\tilde{f}$ are linear maps of equal rank, the
latter being referred to as the {\em curvature rank} at $m$. Let
$B_m$ denote the range space of $f$ or $\tilde{f}$ at $m$ (according
to the agreed identif\/ication) so that $\dim B_m$ equals the
curvature rank at $m$ and which, in turn, is $\leq6$. This leads to
a convenient algebraic classif\/ication of $\rm Riem$ at $m$ into f\/ive
mutually exclusive and disjoint {\em curvature classes} (for further
details, see~\cite{2}).
\begin{enumerate}[{\bf {\it Class}  A}]\itemsep=0pt
\item This covers all possibilities not covered by classes $\mathbf{B}$, $\mathbf{C}$, $\mathbf{D}$ and $\mathbf{O}$ below. For this class, the curvature rank at $m$ is 2, 3, 4, 5 or 6.
\item This occurs when $\dim B_m=2$ and when $B_m$ is spanned by a timelike-spacelike pair of simple bivectors with
orthogonal blades (chosen so that one is the dual of the other). In
this case, one can choose a null tetrad $l,n,x,y\in T_mM$ such that
these bivectors are $F=l\wedge n$ and $\overset{*}{F}=x\wedge y$ so
that $F$ is timelike and $\overset{*}{F}$ is spacelike and then
(using the algebraic identity $R_{a[bcd]}=0$ to remove cross terms)
one has, at $m$,
    \begin{equation}\label{E4}
        R_{abcd}=\alpha F_{ab}F_{cd}+\beta \overset{*}{F}_{ab}\overset{*}{F}_{cd}
    \end{equation}
    for $\alpha,\beta\in\mathbb{R}$, $\alpha\neq0\neq\beta$.
\item In this case $\dim B_m=2$ or $3$ and $B_m$ may be spanned by independent simple bivectors~$F$ and~$G$ (or $F$, $G$ and $H$) with the property that there exists $0\neq r\in T_mM$ such that~$r$ lies in the blades of $\overset{*}{F}$ and $\overset{*}{G}$ (or $\overset{*}{F}$, $\overset{*}{G}$ and $\overset{*}{H}$). Thus $F_{ab}r^b=G_{ab}r^b(=H_{ab}r^b)=0$ and~$r$ is then unique up to a multiplicative non-zero real number.
\item In this case $\dim B_m=1$. If $B_m$ is spanned by the bivector $F$ then, at $m$,
    \begin{equation}\label{E5}
        R_{abcd}=\alpha F_{ab}F_{cd}
    \end{equation}
    for $0\neq\alpha\in\mathbb{R}$ and $R_{a[bcd]}=0$ implies that $F_{a[b}F_{cd]}=0$ from which it may be checked that $F$ is necessarily simple.
\setcounter{enumi}{14}
\item In this case ${\rm Riem}$ vanishes at $m$.
\end{enumerate}
\newpage

It is remarked that this classif\/ication is pointwise and may vary
over $M$. The subset of $M$ consisting of points at which the
curvature class is $\mathbf{A}$ is an open subset of $M$ \cite[p.~393]{2} and the analogous subset arising from the class $\mathbf{O}$
is closed (and has empty interior in the manifold topology of $M$ if
$(M,g)$ is non-f\/lat). It is also useful to note that {\em the equation
$R_{abcd}k^d=0$ at $m$ has no non-trivial solutions if the curvature
class at $m$ is $\mathbf{A}$ or $\mathbf{B}$, a unique independent
solution $($the vector $r$ above$)$ if the curvature class at $m$ is
$\mathbf{C}$ and two independent solutions if the curvature class at
$m$ is $\mathbf{D}$ $($and which span the blade of $\overset{*}{F}$ in
\eqref{E5}}). If $\dim B_m\geq4$ the curvature class at $m$ is
$\mathbf{A}$. If $(M,g)$ has the same curvature class at each $m\in
M$ it will be referred to as being \emph{of that class}.

Related to this classif\/ication scheme is the following result which will prove useful in what is to follow. The details and proof can be found in \cite{3,4,5,2}.
\begin{theorem}\label{theorem1}
Let $(M,g)$ be a space-time, let $m\in M$ and let $h$ be a second
order, symmetric, type $(0,2)$ $($not necessarily non-degenerate$)$ tensor
at $m$ satisfying $h_{ae}R^e_{\ bcd}+h_{be}R^e_{\ acd}=0$. Then,
with all tensor index movements and all orthogonality statements
made using the metric $g$;
\begin{enumerate}[$(i)$]
\item if the curvature class of $(M,g)$ at $m\in M$ is $\mathbf{D}$ and $u,v\in T_mM$ span the $2$-space at $m$ orthogonal to $F$ in \eqref{E5} $($that is $u\wedge v$ is the blade of $\overset{*}{F})$ there exists $\phi,\mu,\nu,\lambda\in\mathbb{R}$ such that, at
$m$,
    \begin{equation}\label{E6}
        h_{ab}=\phi g_{ab}+\mu u_au_b+\nu v_av_b+\lambda(u_av_b+v_au_b);
    \end{equation}
\item if the curvature class of $(M,g)$ at $m\in M$ is $\mathbf{C}$ there exists $r\in T_mM$ $($the vector appearing in the above definition of class $\mathbf{C})$ and $\phi,\lambda\in\mathbb{R}$ such that, at
$m$,
    \begin{equation}\label{E7}
        h_{ab}=\phi g_{ab}+\lambda r_ar_b;
    \end{equation}
\item if the curvature class of $(M,g)$ at $m\in M$ is $\mathbf{B}$ there exists a null tetrad $l$, $n$, $x$, $y$ $($that appearing in the above definition of class $\mathbf{B})$ and $\phi,\lambda\in\mathbb{R}$ such that, at
$m$,
    \begin{equation}\label{E8}
        h_{ab}=\phi g_{ab}+\lambda(l_an_b+n_al_b)=(\phi+\lambda)g_{ab}-\lambda(x_ax_b+y_ay_b);
    \end{equation}
\item if the curvature class of $(M,g)$ at $m\in M$ is $\mathbf{A}$ there exists $\phi\in\mathbb{R}$ such that, at
$m$,
    \begin{equation}\label{E9}
        h_{ab}=\phi g_{ab}.
    \end{equation}
\end{enumerate}
\end{theorem}

The proof is essentially based on the obvious fact that the range
$B_m$ of the map $f$ in (\ref{E3}) which, by the algebraic
symmetries of the curvature $\rm Riem$ of $(M,g)$, consists entirely of
members which are skew-self adjoint with respect to $g$, must
likewise consist entirely of members skew-self adjoint with respect
to $h$. Thus each $F\in B_m$ satisf\/ies
\begin{equation}\label{E10}
g_{ac}F^c_{\ b}+g_{bc}F^c_{\ a}=0, \qquad h_{ac}F^c_{\
b}+h_{bc}F^c_{\ a}=0.
\end{equation}
It is a consequence of (\ref{E10}) that the blade of $F$ (if $F$ is
simple) and each of the canonical pair of blades of $F$ (if $F$ is
non-simple) are eigenspaces of $h$ with respect to $g$, that is, for
$F$ simple, any $k$ in the blade of $F$ satisf\/ies $h_{ab}k^b=\omega
g_{ab}k^b$ where the eigenvalue $\omega\in\mathbb{R}$ is independent
of~$k$, and similarly for each of the canonical blades if $F$ is
non-simple (but with possibly dif\/ferent eigenvalues for these
blades) \cite{2,3,4}.

It is remarked that if $g'$ is another metric on $M$ whose curvature
tensor ${\rm Riem}'$ equals the curvature tensor ${\rm Riem}$ of $g$
\emph{everywhere on $M$} then the conditions of this theorem are
satisf\/ied for $h=g'(m)$ at each $m\in M$ and so the conclusions also
hold except that now one must add the restriction $\phi\neq 0$ in
each case to preserve the non-degeneracy of $g'$ at $m$ and maybe
some restrictions on $\phi,\mu,\nu$ and $\lambda$ if the signature
of $g'$ is prescribed. If $(M,g)$ is of class $\mathbf{A}$,
(\ref{E9})~gives $g'=\phi g$ and the Bianchi identity may be used to
show that $\phi$ is constant on $M$~\cite{4,2}. When Theorem~\ref{theorem1} is
applied to another metric $g'$ on $M$ in this way, it consolidates
the curvature classif\/ication scheme which preceded it. To see this
note that $g'$ need not have Lorentz signature $(-,+,+,+)$ but, if
this is insisted upon, the curvature classif\/ication scheme,
including the nature (timelike, spacelike or null) of $F$ and
$\overset{*}{F}$ in class $\mathbf{B}$, $r$ in class $\mathbf{C}$
and $F$ in class $\mathbf{D}$ is the same whether taken for (the
common curvature tensor) $\rm Riem$ with $g$ or $\rm Riem$ with $g'$. In
this sense it is a classif\/ication of $\rm Riem$, independent of the
metric generating $\rm Riem$ \cite{6}.

\section{Holonomy theory}\label{section4}

Let $(M,g)$ be a space-time with Levi-Civita connection $\nabla$ and
let $m\in M$. The connection $\nabla$ is a complicated object but
one particularly pleasing feature of it lies in the following
construction. Let $m\in M$ and for $1\leq k\leq\infty$ let $C_k(m)$
denote the set of all piecewise $C^k$ closed curves starting and
ending at $m$. If $c\in C_k(m)$ let $\tau_c$ denote the vector space
isomorphism of $T_mM$ obtained by parallel transporting, using
$\nabla$, each member of $T_mM$ along $c$. Using a standard notation
associated with curves one def\/ines, for curves $c,c_0,c_1,c_2\in
C_k(m)$, with $c_0$ denoting a constant curve at $m$, the {\em
identity map} $\tau_{c_0}$ on $T_mM$, the {\em inverse}
$\tau^{-1}_c\equiv\tau_{c^{-1}}$ and {\em product}
$\tau_{c_1}\cdot\tau_{c_2}\equiv\tau_{c_1\cdot c_2}$ to put a group
structure on $\{\tau_c:c\in C_k(m)\}$, making it a subgroup of $G\equiv
GL(T_mM)(=GL(4,\mathbb{R}))$, called the {\em $k$-holonomy group} of
$M$ at $m$ and denoted by $\Phi_k(m)$. In fact, since $M$ is
connected and also a manifold, it is also path connected and, as a
consequence, it is easily checked that, up to an isomorphism,
$\Phi_k(m)$ is independent of $m$. Less obvious is the fact that
$\Phi_k(m)$ is independent of $k$ ($1\leq k\leq\infty$) and thus one
arrives at the {\em holonomy group} $\Phi$ (of $\nabla$) on $M$.
Further details may be found in \cite{7} and a summary in \cite{2}.

One could repeat the above operations, but now only using curves
homotopic to zero. the above independence of $m$ and $k$ still holds
and one arrives at the {\em restricted holonomy group} $\Phi^0$ of
$M$. It can now be proved that $\Phi$ and $\Phi^0$ are Lie subgroups
of $G$, with $\Phi^0$ connected, and that $\Phi^0$ is the identity
component of $\Phi$ \cite{7}. Clearly, if $M$ is simply connected,
$\Phi=\Phi^0$ and then $\Phi$ is a connected Lie subgroup of $G$.
The common Lie algebra $\phi$ of $\Phi$ and $\Phi^0$ is called the
{\em holonomy algebra}. The connection $\nabla$ can then be shown to
be {\em flat} (that is, $\rm Riem$ vanishes on $M$) if and only if
$\phi$ is trivial and this, in turn, is equivalent to $\Phi^0$ being
trivial.

For a space-time $(M,g)$, however, one has the additional
information that $\nabla$ is compatible with the metric $g$, that
is, $\nabla g=0$. Thus each map $\tau_c$ on $T_mM$ preserves inner
products with respect to $g(m)$. It follows that $\Phi$ is
(isomorphic to) a subgroup of the {\em Lorentz group} $\mathcal{L}$,
where $\mathcal{L}=\{A\in GL(4,\mathbb{R}):A\eta A^T=\eta\}$ where
$A^T$ denotes the transpose of $A$ and $\eta$ is the Minkowski
metric, $\eta=\mathrm {diag}(-1,1,1,1)$. Now $\Phi$ is a Lie
subgroup of $GL(4,\mathbb{R})$ and $\mathcal{L}$ can be shown to be
a 6-dimensional Lie subgroup of $GL(4,\mathbb{R})$ which is a {\em
regular} submanifold of $GL(4,\mathbb{R})$ (that is, the
(sub)manifold topology on $\mathcal{L}$ equals its induced topology
from $GL(4,\mathbb{R})$). It follows that $\Phi$ is a~{\em Lie
subgroup} of $\mathcal{L}$ (see e.g.~\cite{2}) and hence that
$\Phi^0$ (or $\Phi$, if $M$ is simply connected) is a {\em
connected} Lie subgroup of the identity component, $\mathcal{L}_0$,
of $\mathcal{L}$. Thus the holonomy algebra $\phi$ can be identif\/ied
with a subalgebra of the Lie algebra $L$ of $\mathcal{L}$, the {\em
Lorentz algebra}. The one-to-one correspondence between the
subalgebras of $L$ and the connected Lie subgroups of
$\mathcal{L}_0$ shows that the Lie group~$\Phi^0$ (or $\Phi$, if $M$
is simply connected) is determined by the subalgebra of $L$
associated with $\phi$. Fortunately, the subalgebra structure of $L$
is well-known and can be conveniently represented as follows. Let
$m\in M$ and choose a basis for $T_mM$ (together with its dual
basis) for which $g(m)$ has components equal to, say, $g_{ab}$. Then
$\mathcal{L}$ is isomorphic to $\{A\in GL(4,\mathbb{R}):Ag(m)
A^T=g(m)\}$ and $L$ can then be represented as the subset of
$M_4\mathbb{R}$ each member of which, if regarded as the set of
components of a type $(1,1)$ tensor at $m$ in this basis, is
skew-self adjoint with respect to the matrix $g_{ab}(m)$ (that is,
its components $F^a_{\ b}$ satisfy an equation like the f\/irst in~(\ref{E10}) with respect to $g(m)$). Thus one can informally
identify~$L$ with this collection of bivectors, usually, written in
type $(2,0)$ form. The binary operation on~$L$ is that induced from
the Lie algebra $M_4\mathbb{R}$ of $GL(4,\mathbb{R})$ and is matrix
commutation. Such a~representation of $L$ is well-known and has been
classif\/ied into f\/ifteen convenient types~\cite{8} (for details of
the possible holonomy types most relevant for the physics of general
relativity see~\cite{25}). It is given in the f\/irst three columns of
Table~\ref{table1} using either a null tetrad $l$, $n$, $x$, $y$ or an orthonormal
tetrad $u$, $x$, $y$, $z$ to describe a basis for each subalgebra.
\begin{table}\centering
\caption{Holonomy algebras. For types $R_5$ and $R_{12}$
$0\neq\omega\in\mathbb{R}$.
Every potential holonomy algebra except $R_5$ (for which the curvature tensor fails to satisfy the algebraic Bianchi identity) can occur as an actual holonomy algebra,
see, e.g.~\cite{2}). There is also a type $R_1$, when $\phi$ is trivial and
$(M,g)$ is f\/lat, but this trivial type is omitted. Type $R_{15}$ is
the ``general'' type, when $\phi=L$.}\label{table1}
\vspace{1mm}

\begin{tabular}{|@{\,\,}c@{\,\,}|@{\,\,}c@{\,\,}|@{\,\,}c@{\,\,}|@{\,\,}c@{\,\,}|@{\,\,}c@{\,\,}|@{\,\,}c@{\,\,}|}
\hline
Type & Dimension & Basis & Curvature & Recurrent & Constant \cr
     &           &       & Class     & vector f\/ields & vector f\/ields \cr
\hline
 \tsep{0.5mm}$R_{2}$ & 1 & $l\wedge n$ & $\mathbf{D}$ or $\mathbf{O}$ &$\{l\}$, $\{n\}$ & $\langle x, y\rangle $ \\
 $R_{3}$ & 1 & $l\wedge x$ & $\mathbf{D}$ or $\mathbf{O}$ & -- & $\langle l, y\rangle $ \\
 $R_{4}$ & 1 & $x\wedge y$ & $\mathbf{D}$ or $\mathbf{O}$ &-- & $\langle l, n\rangle$ \\
 $R_{5}$ & 1 & $l\wedge n+\omega x\wedge y$ & -- &-- & -- \\
 $R_{6}$ & 2 & $l\wedge n$, $l\wedge x$ & $\mathbf{C}$, $\mathbf{D}$ or $\mathbf{O}$ & $\{l\}$ & $\langle y\rangle$ \\
 $R_{7}$ & 2 & $l\wedge n$, $x\wedge y$ & $\mathbf{B}$, $\mathbf{D}$ or $\mathbf{O}$ & $\{l\}$, $\{n\}$ & -- \\
 $R_{8}$ & 2 & $l\wedge x$, $l\wedge y$ & $\mathbf{C}$, $\mathbf{D}$ or $\mathbf{O}$ & -- & $\langle l\rangle$ \\
 $R_{9}$ & 3 & $l\wedge n$, $l\wedge x$, $l\wedge y$ & $\mathbf{A}$, $\mathbf{C}$, $\mathbf{D}$ or $\mathbf{O}$ & $\{l\}$ & -- \\
 $R_{10}$ & 3 & $l\wedge n$, $l\wedge x$, $n\wedge x$ & $\mathbf{C}$, $\mathbf{D}$ or $\mathbf{O}$ & -- & $\langle y\rangle$ \\
 $R_{11}$ & 3 & $l\wedge x$, $l\wedge y$, $x\wedge y$ & $\mathbf{C}$, $\mathbf{D}$ or $\mathbf{O}$ & -- & $\langle l\rangle$ \\
 $R_{12}$ & 3 & $l\wedge x$, $l\wedge y$, $l\wedge n+\omega(x\wedge y)$ & $\mathbf{A}$, $\mathbf{C}$, $\mathbf{D}$ or $\mathbf{O}$ & $\{l\}$ & -- \\
 $R_{13}$ & 3 & $x\wedge y$, $y\wedge z$, $x\wedge z$ & $\mathbf{C}$, $\mathbf{D}$ or $\mathbf{O}$ & -- & $\langle u\rangle$ \\
 $R_{14}$ & 4 & $l\wedge n$, $l\wedge x$, $l\wedge y$, $x\wedge y$ & any & $\{l\}$ & -- \\
 $R_{15}$ & 6 & $L$ & any &-- & --  \\
 \hline
\end{tabular}
\end{table}

Now suppose $M$ is simply connected. This condition is not always
required but is imposed in this section for convenience. It can
always be assumed in local work, for example, in some connected,
simply connected coordinate domain. In this case $\Phi=\Phi^0$ and
is connected. In any case $\Phi$ will be referred to according to its Lie
algebra label as in Table~\ref{table1}. Two important results can now be
mentioned in connection with the f\/irst column of Table~\ref{table1}. First, it
turns out \cite{9,10,2} that if $m\in M$ there exists $0\neq k\in
T_mM$ such that $F^a_{\ b}k^b=0$ for each $F\in\phi$ if and only if
$M$ admits a global, {\em covariantly constant}, smooth vector f\/ield
whose value at $m$ is $k$. A basis for the vector space of such
vector f\/ields on $M$ for each holonomy type is given inside $\langle \ \rangle$
brackets in the f\/inal column of Table~\ref{table1}. Second, there exists $0\neq
k\in T_mM$ such that $k$ is an eigenvector of each $F\in\phi$ but
with at least one associated eigenvalue not zero if and only if $M$
admits a global smooth {\em properly recurrent} vector f\/ield $X$
whose value at $m$ is $k$, that is, a global nowhere zero vector
f\/iel~ $X$ on~$M$ satisfying $\nabla X=X\otimes w$ for some global,
smooth covector f\/ield $w$ on~$M$ (the {\em recurrence $1$-form}) and
such that no function $\alpha:M\rightarrow\mathbb{R}$ exists such
that $\alpha$ is nowhere zero on~$M$ and $\alpha X$ is covariantly
constant on~$M$. It follows from the existence of one non-zero
eigenvalue at $m$ in the above def\/inition that $X$ is necessarily
{\em null} (by the skew-self adjoint property of the members of
$\phi$). The independent {\em properly recurrent} vector f\/ields are
listed for each holonomy type in $\{\ \}$ brackets in the second
from last column of Table~\ref{table1}. (It is remarked here that a nowhere
zero vector f\/ield $Y$ on $M$ is called {\em recurrent} if $\nabla
Y=Y\otimes r$ for some global covector f\/ield $r$ on $M$. In this
case $r$ could be identically zero and so (non-trivial) covariantly
constant vector f\/ields are, in this sense, recurrent. In fact, any
non-null recurrent vector f\/ield or any recurrent vector f\/ield on a
manifold with positive def\/inite metric can be globally scaled to be
nowhere zero and covariantly constant because if $Y$ is any such
vector f\/ield and $\nabla Y=Y\otimes r$, $\alpha Y$ is covariantly
constant, where $\alpha=\exp(-\tfrac{1}{2}\log |g(Y,Y)|)$.)

A recurrent vector f\/ield is easily seen to def\/ine a 1-dimensional
distribution on $M$ which is preserved by parallel transport. There
is an important generalisation of this concept. Let $m\in M$ and $V$
a non-trivial proper subspace of $T_mM$. Suppose $\tau_c(V)=V$ for
each $\tau_c$ arising from $c\in C_k(m)$ at $m$. Then $V$ is {\em
holonomy invariant} and gives rise in an obvious way to a~smooth
distribution on $M$ which is, in fact, integrable~\cite{7}. Clearly,
if $V\subset T_mM$ is holonomy invariant then so is the orthogonal
complement, $V^{\perp}$, of~$V$. If such a $V$ exists the holonomy
group $\Phi$ of $M$ is called {\em reducible} (otherwise, {\em
irreducible}). (This concept of holonomy reducibility is a little
more complicated in the case of a Lorentz metric than in the
positive def\/inite case due to the possibility of null holonomy
invariant subspaces giving a weaker form of reducibility. This will
not be pursued any further here, more details being available in
\cite{11} and summaries in \cite{2,25,10}.) Thus, for example, in
the notation of Table~\ref{table1} the holonomy type $R_2$ admits two
1-dimensional null holonomy invariant subspaces spanned by $l$ and
$n$ and which give rise to two null properly recurrent vector f\/ields
and inf\/initely many 1-dimensional spacelike holonomy invariant
subspaces spanned by the inf\/initely many covariantly constant vector
f\/ields in $\langle x,y\rangle$. For the holonomy type $R_7$, two 1-dimensional
null holonomy invariant subspaces exist and which give rise to two
independent properly recurrent null vector f\/ields as in the previous
case, together with a 2-dimensional spacelike one orthogonal to each
of the null ones. For holonomy types $R_{10}$, $R_{11}$ and $R_{13}$
one has a 1-dimensional holonomy invariant subspace, spanned by a
covariantly constant vector f\/ield in each case, together with its
orthogonal complement. These holonomy decompositions will be useful
in Sections~\ref{section6} and~\ref{section7}. In the event that $M$ is not simply
connected, the vector f\/ields determined by the holonomy and
described above may not exist globally but do exist locally over some open, connected and simply
connected neighbourhood of any point.

It is remarked here that if $(M,g)$
is of curvature class $\mathbf{B}$ it can be shown that it must in fact be
of holonomy type $R_7$ \cite{12,2}.

It is useful, at this point, to introduce the {\em infinitesimal
holonomy group} $\Phi'_m$, of $(M,g)$ at each $m\in M$. Using a
semi-colon to denote a $\nabla$-covariant derivative, consider, in
some coordinate neighbourhood of $m$, the following matrices for
$(M,g)$ at $m$
\begin{equation}\label{E11}
R^a_{\ bcd}X^cY^d,\qquad R^a_{\ bcd;e}X^cY^dZ^e,\qquad \ldots
\end{equation}
for $X,Y,Z,\ldots \in T_mM$. It turns out that the collection
(\ref{E11}) spans a subalgebra of the holonomy algebra $\phi$ (and
hence only a f\/inite number of terms arise in (\ref{E11})~\cite{7}).
This algebra is called the {\em infinitesimal holonomy algebra at
$m$} and is denoted by $\phi'_m$. The unique connected Lie subgroup
of $\Phi$ that it gives rise to is the {\em infinitesimal holonomy
group} $\Phi'_{m}$ at $m$. This is useful in that it says that the
range space of the map $f$ in (\ref{E3}) is, at each $m\in M$,
(isomorphic as a vector space to) a subspace of $\phi$. This gives a
restriction, when $\phi$ is known, on the expression for $\rm Riem$ at
each $m$ and hence on its curvature class at $m$. This restriction
is listed in the fourth column of Table~\ref{table1}. Thus if the holonomy type
of $(M,g)$ is $R_2$, $R_3$ or $R_4$ its curvature class is $\mathbf{O}$ or
$\mathbf{D}$ at each $m\in M$ whilst if it is $R_6$, $R_8$, $R_{10}$,
$R_{11}$ or $R_{13}$ it is $\mathbf{O}$, $\mathbf{D}$ or $\mathbf{C}$, if $R_7$ it is $\mathbf{O}$, $\mathbf{D}$
or $\mathbf{B}$, if $R_9$ or $R_{12}$ it is $\mathbf{O}$, $\mathbf{D}$, $\mathbf{C}$ or $\mathbf{A}$ and if
$R_{14}$ or~$R_{15}$ it could be any curvature class. A useful
relationship between the various algebras $\phi'_m$, the algebra
$\phi$ and the curvature class (through the range space~ $B_{m}$) at
each $m\in M$ is provided by the Ambrose--Singer theorem \cite{13}
(see also~\cite{7}).

\section{Projective structure}\label{section5}

One aspect of dif\/ferential geometry that has been found interesting
both for pure geometers and physicists working in general relativity
theory is that of projective structure. For general relativity it is
clearly motivated by the Newton--Einstein principle of equivalence.
In this section it will take the form described in the following
question; for a space-time $(M,g)$ with Levi-Civita connection
$\nabla$, if one knows the paths of all the {\em unparameterised}
geodesics (that is, only the {\em geodesic paths} in $M$) how
tightly is $\nabla$ determined? Put another way, let $(M,g)$ and
$(M,g')$ be space-times with respective Levi-Civita connections
$\nabla$ and $\nabla'$, such that the sets of geodesic paths of
$\nabla$ and $\nabla'$ coincide (and let it be agreed that $\nabla$
and $\nabla'$ (or $g$ and $g'$) are then said to be {\em
projectively related} (on $M$)). What can be deduced about the
relationship between $\nabla$ and $\nabla'$ (and between $g$ and
$g'$)? (And it is, perhaps, not surprising that some reasonable link
should exist between projective relatedness and holonomy theory.) In
general, $\nabla$ and $\nabla'$ may be expected to dif\/fer but it
turns out that in many interesting situations they are necessarily
equal. If $\nabla=\nabla'$ is the result, holonomy theory can also
describe precisely, the (simple) relationship between $g$ and $g'$~\cite{14,2}.

If $\nabla$ and $\nabla'$ are projectively related then there exists
a uniquely def\/ined global smooth 1-form f\/ield $\psi$ on $M$ such
that, in any coordinate domain of $M$, the respective Christof\/fel
symbols of $\nabla$ and $\nabla'$ satisfy \cite{15,16}
\begin{equation}\label{E12}
\Gamma'^a_{\ bc}-\Gamma^a_{bc}=\delta^a_{\ b}\psi_c+\delta^a_{\ c}\psi_b.
\end{equation}
It is a consequence of the fact that $\nabla$ and $\nabla'$ are {\em
metric connections} that $\psi$ is a {\em global} gradient on $M$
\cite{15}. Equation (\ref{E12}) can, by using the identity
$\nabla'g'=0$, be written in the equivalent form
\begin{gather}\label{E13}
g'_{ab;c}=2g'_{ab}\psi_c+g'_{ac}\psi_b+g'_{bc}\psi_a
\end{gather}
(recalling that a semi-colon denotes covariant dif\/ferentiation with respect to $\nabla$).
Equation (\ref{E12}) reveals a simple relation between the type $(1,3)$ curvature tensors ${\rm Riem}$ and ${\rm Riem}'$ of~$\nabla$ and~$\nabla'$, respectively, given by
\begin{equation}\label{E14}
R'^a_{\ bcd}=R^a_{\ bcd}+\delta^a_{\ d}\psi_{bc}-\delta^a_{\ c}\psi_{bd} \hspace{1cm} (\Rightarrow R'_{ab}=R_{ab}-3\psi_{ab}),
\end{equation}
where $\psi_{ab}\equiv \psi_{a;b}-\psi_a\psi_b=\psi_{ba}$ and where
$R'_{ab}\equiv R'^c_{\ acb}$ are the Ricci tensor components of
$\nabla'$. It can now be shown that if $\nabla$ and $\nabla'$ are
projectively related, the following type $(1,3)$ {\em Weyl
projective tensor}, $W$, is the same for each of them~\cite{17}
\begin{equation}\label{E15}
W^a_{\ bcd}=R^a_{\ bcd}+\tfrac{1}{3}(\delta^a_{\ d}R_{bc}-\delta^a_{\ c}R_{bd}).
\end{equation}

A particularly important case of such a study arises where the
original pair $(M,g)$ is a~space-time which is also an {\em Einstein
space} so that the tensor $E$ in (\ref{E2}) is identically zero on~$M$. Such a situation has been discussed in several places \cite{18,19,20,21,22,23}  
in connection with the principle of equivalence. The particular case
which is, perhaps, of most importance in general relativity arises
when the Ricci scalar vanishes and then $(M,g)$ is a vacuum (Ricci
f\/lat) space-time and this is discussed in \cite{21, 23}. It turns
out that if $(M,g)$ is a space-time which is a~(general) Einstein
space (and with the non-f\/lat assumption temporarily dropped) and if~$g'$ is another metric on~$M$ projectively related to~$g$ then
either $(M,g)$ and $(M,g')$ are each of constant curvature, or the
Levi-Civita connections $\nabla$ and $\nabla'$ of~$g$ and~$g'$,
respectively, are equal. In the event that neither space-time is of
constant curvature, (and so $\nabla'=\nabla$) an argument from
holonomy theory can be used to show that, \emph{generically},
$(M,g')$ is also an Einstein space and that $g'=cg$ ($0\neq
c\in\mathbb{R}$). Although examples exist where each of these
conclusions fail, $g'$~always has Lorentz signature (up to an
overall minus sign). If, in addition, $(M,g)$ is assumed vacuum and
the non-f\/lat condition is imposed, then, necessarily,
$\nabla=\nabla'$ and, with one very special case excluded, $g'=cg$
on $M$ ($0\neq c\in\mathbb{R}$) (and so $(M,g')$ is also vacuum).
For this case $g'$ has the same signature as $g$ (up to an overall
minus sign \cite{21, 23}). This result is relevant for general
relativity theory. A similar restrictive result for space-times of
certain holonomy types will be established in the next two sections.

The formalism described above characterises the relationship between
the connections and metrics of two projectively related space-times
$(M,g)$ and $(M,g')$. However, it is convenient to rewrite them in a
dif\/ferent way using the Sinyukov transformation \cite{19}. Thus,
with this projective relatedness assumed, one takes advantage of the
fact that the 1-form $\psi$ in (\ref{E12}) is necessarily a global
gradient by writing $\psi=d\chi$ for some smooth function
$\chi:M\rightarrow\mathbb{R}$. Then the pair $g'$ and $\psi$ above
are replaced by a type $(0,2)$ symmetric tensor f\/ield $a$ and a
1-form f\/ield~$\lambda$ on~$M$ which are given in terms of $g'$ and
$\psi$ by
\begin{equation}\label{E16}
a_{ab}=e^{2\chi}g'^{cd}g_{ac}g_{bd},\qquad \lambda_a=-e^{2\chi}\psi_bg'^{bc}g_{ac} \quad (\Rightarrow \lambda_a=-a_{ab}\psi^b),
\end{equation}
where an abuse of notation has been used in that $g'^{ab}$ denotes
the contravariant components of~$g'$ (and not the tensor $g'_{ab}$
with indices raised using $g$) so that $g'_{ac}g'^{cb}=\delta_a^{\
b}$. Then (\ref{E16}) may be inverted to give
\begin{equation}\label{E17}
g'^{ab}=e^{-2\chi}a_{cd}g^{ac}g^{bd},\qquad
\psi_a=-e^{-2\chi}\lambda_bg^{bc}g'_{ac}.
\end{equation}
The idea is that if $g$ and $g'$ are projectively related metrics on
$M$, so that (\ref{E13}) holds for some 1-form $\psi(=d\chi)$, then
$a$ and $\lambda$ as def\/ined in (\ref{E16}) can be shown, after a
short calculation, to satisfy Sinjukov's equation
\begin{equation}\label{E18}
a_{ab;c}=g_{ac}\lambda_b+g_{bc}\lambda_a.
\end{equation}
From (\ref{E16}) it follows that $a$ is non-degenerate and  from
(\ref{E18}), after a contraction with $g^{ab}$, that $\lambda$ is a
global gradient on $M$ (in fact, of $\tfrac{1}{2}a_{ab}g^{ab}$).

In practice, when asking which pairs $(g',\nabla')$ are projectively
related to some original pair $(g,\nabla)$ on $M$, it is often
easier to use (\ref{E18}) to attempt to f\/ind $a$ and $\lambda$
rather than (\ref{E13}) to f\/ind $g'$ and $\psi$. But then one must
be able to convert back from $a$ and $\lambda$ to $g'$ and $\psi$.
To do this one f\/irst assumes that such a non-degenerate tensor $a$
and a 1-form $\lambda$ are given on $M$ and which together satisfy
(\ref{E18}) and so one may def\/ine a symmetric non-degenerate type
$(2,0)$ tensor $a^{-1}$ on $M$ which, at each $m\in M$, is the
inverse of $a$ ($a_{ac}(a^{-1})^{cb}=\delta_a^{\ b}$). Then raising
and lowering indices on $a$ and $a^{-1}$ with $g$ in the usual way,
so that $a^{-1}_{ac}a^{cb}=\delta_a^{\ b}$, one def\/ines a global
1-form $\psi$ on $M$ by $\psi_a=-a^{-1}_{ab}\lambda^b$ (and so
$\lambda_a=-a_{ab}\psi^b$). It follows that \emph{$\psi$ is a global
gradient on $M$}. To see this f\/irst dif\/ferentiate the condition
$a^{-1}_{ac}a^{cb}=\delta_a^{\ b}$ and use (\ref{E18}) to f\/ind,
after a short calculation,
\begin{equation}\label{E19}
a^{-1}_{ab;c}=a^{-1}_{ac}\psi_b+a^{-1}_{bc}\psi_a.
\end{equation}
Now def\/ine a smooth, symmetric connection $\nabla''$ on $M$ by
decreeing that, in any coordinate domain, its Christof\/fel symbols
are given by $\Gamma''$ where
\begin{equation}\label{E20}
\Gamma''^a_{\ bc}=\Gamma^a_{bc}-\psi^ag_{bc}.
\end{equation}
Then it is easily checked from (\ref{E19}), (\ref{E20}) and the
equation $\lambda_a=-a_{ab}\psi^b$ that $\nabla''a=0$ and so if the
everywhere non-degenerate tensor $a$ is regarded as a metric on $M$,
$\nabla''$ is its Levi-Civita connection. But a contraction of
(\ref{E20}) over the indices $a$ and $c$ gives, from a standard
formula for Christof\/fel symbols
\begin{equation}\label{E21}
\Gamma''^a_{\ ba}-\Gamma^a_{ba}=\tfrac{\partial}{\partial x^b}\left(
\tfrac{1}{2}\ln\left(\frac{|\det a|}{|\det g|}\right)\right)=-\psi_b
\end{equation}
and so $\psi$ is seen to be a \emph{global gradient} on $M$. (The authors have recently discovered that Sinyukov \cite{19} had established the same result by a method involving the direct construction of the connection of the tensor $a^{-1}$.) Writing
$\chi$ for this potential function, so that $\psi=d\chi$, one
def\/ines a metric $g'$ on $M$ by $g'=e^{2\chi}a^{-1}$. Then (\ref{E19})
can be used to show that $g'$ and $\psi$ satisfy~(\ref{E13}) and
hence that $g'$ is projectively related to $g$. It is easily checked
that the tensors~$\psi$ and~$g'$ thus found satisfy (\ref{E16}) and
(\ref{E17}). [It is remarked that if $\psi$ is replaced by $\lambda$
and the metric~$g$ by the tensor $a^{-1}$ in (\ref{E20}) the
connection $\nabla'''$ thus def\/ined satisf\/ies $\nabla'''a^{-1}=0$.
In fact, $\nabla=\nabla'' \Leftrightarrow \nabla=\nabla'''
\Leftrightarrow \nabla=\nabla'$.] Since any solution pair
$(g',\psi)$ of (\ref{E13}) leads to a~pair $(a,\lambda)$ satisfying
(\ref{E18}) it follows that \emph{all projectively related metrics~$g'$ together with their attendant $1$-forms $\psi$ will be found if
all pairs $(a,\lambda)$ can be found and which, together, satisfy~\eqref{E18}}.

For the f\/inding of the general solution of (\ref{E18}) a useful
result arises by applying the Ricci identity to $a$ and using
(\ref{E18}) to get
\begin{equation}\label{E22}
(a_{ab;cd}-a_{ab:dc}=)a_{ae}R^e_{\ bcd}+a_{be}R^e_{\ acd}=g_{ac}\lambda_{bd}+g_{bc}\lambda_{ad}-g_{ad}\lambda_{bc}-g_{bd}\lambda_{ac},
\end{equation}
where $\lambda_{ab}=\lambda_{a;b}=\lambda_{ba}$. This leads to the
following lemma (which is a special case of a more detailed result
in \cite{21}) and for which a def\/inition is required. Suppose $m\in
M$, that $F\in\Lambda_mM$ and that the curvature tensor ${\rm Riem}$ of
$(M,g)$ satisf\/ies $R^{ab}_{\ \ cd}F^{cd}=\alpha F^{ab}$
($\alpha\in\mathbb{R}$) at $m$ so that $F$ is a (real) eigenvector
of the map $\tilde{f}$ in (\ref{E3}). Then $F$ is called a (real)
{\em eigenbivector of ${\rm Riem}$ at~$m$ with eigenvalue $\alpha$}.

\begin{lemma}\label{lemma1}
Let $(M,g)$ and $(M,g')$ be space-times with $g$ and $g'$
projectively  related. Suppose at $m\in M$ that $F\in\Lambda_mM$ is
a $($real$)$ eigenbivector of ${\rm Riem}$ of $(M,g)$ with zero eigenvalue $($so
that $F$ is in the kernel, $\ker\tilde{f}$, of $\tilde{f}$ in
\eqref{E3}$)$. Then the blade of $F$ (if $F$ is simple) or each of the
canonical pair of blades of $F$ $($if $F$ is non-simple$)$ is an
eigenspace of the symmetric tensor $\lambda_{ab}$ with respect to
$g$ at $m$. $($That is, if $p\wedge q$ is in any of these blades
$(p,q\in T_mM)$ there exists $\mu\in\mathbb{R}$ such that for any
$k\in p\wedge q$, $\lambda_{ab}k^b=\mu g_{ab}k^b)$. In particular,
\begin{enumerate}[$(i)$]\itemsep=0pt
\item suppose that, at $m$,  the collection of all the blades of all simple members of $\ker\tilde{f}$ and all the canonical blade pairs of all the non-simple members of $\ker\tilde{f}$ $($and which are each eigenspaces of $\lambda_{a;b})$ are such that they force $T_mM$ to be an eigenspace of $\lambda_{a;b}$. Then  $\lambda_{a;b}$ is proportional to $g_{ab}$ at $m$.

\item If condition $(i)$ is satisfied at
    each $m'$ in some connected open neighbourhood $U$ of $m$, then, on $U$ and for some $c\in\mathbb{R}$,
    \[ (a)\ \lambda_{ab}=cg_{ab}, \qquad (b)\ \lambda_dR^d_{\ abc}=0, \qquad (c)\ a_{ae}R^e_{\ bcd}+a_{be}R^e_{\ acd}=0.
     \]
\end{enumerate}
\end{lemma}

\begin{proof}
First contract (\ref{E22}) with $F^{cd}$ and use $R^a_{\ bcd}F^{cd}=0$ to get
\begin{equation}\label{E23}
g_{ae}F^{ce}\lambda_{bc}+g_{be}F^{ce}\lambda_{ac}-g_{ae}F^{ed}\lambda_{bd}-g_{be}F^{ed}\lambda_{ad}=0,
\end{equation}
which rearranges, after cancellation of a factor 2, as
\begin{equation}\label{E24}
\lambda_{ac}F^c_{\ b}+\lambda_{bc}F^c_{\ a}=0,
\end{equation}
where $F^a_{\ b}=F^{ac}g_{cb}$ is skew-self adjoint with respect to
$g$.  The argument following (\ref{E10}) applied to
$h_{ab}=\lambda_{ab}$ then completes the proof of the f\/irst part of
the lemma. Part $(i)$ of the lemma then also follows if $T_mM$ is an
eigenspace of $\lambda_{ab}$ with respect to $g$. [It is remarked
here that, for such a symmetric tensor, the eigenspaces
corresponding to distinct eigenvalues are orthogonal.] For part
$(ii)$ one notes that, from part $(i)$, $\lambda_{ab}=\sigma g_{ab}$
on $U$ for some smooth function $\sigma:U\rightarrow\mathbb{R}$. The
Ricci identity on $\lambda$ then shows that
\begin{equation}\label{E25}
\lambda_dR^d_{\ abc}=\lambda_{a;bc}-\lambda_{a;cb}=g_{ab}\sigma_{,c}-g_{ac}\sigma_{,b},
\end{equation}
where a comma denotes a partial derivative. Now for $m'\in U$ there
exists  $0\neq F\in\ker\tilde{f}$ at $m'$ and a contraction of
(\ref{E25}) with  $F^{bc}$ gives
\begin{equation}\label{E26}
g_{ab}(F^{bc}\sigma_{,c})-g_{ac}(F^{bc}\sigma_{,b})=0
\end{equation}
from which it follows that $g_{ab}F^{bc}\sigma_{,c}=0$ and so
$F^{ab}\sigma_{,b}=0$. So either $F$ is non-simple (and hence
$d\sigma=0$) at $m'$ or each member of $\ker\tilde{f}$ is simple and
the 1-form $d\sigma$ is $g$-orthogonal to its blade. In the latter
case, if dim($\ker\tilde{f})=1$, and $\ker\tilde{f}$ is spanned by a
single simple bivector the conditions of $(i)$ are not satisf\/ied,
whereas, if dim($\ker\tilde{f})=2$, $\ker\tilde{f}$ must be spanned
by two simple bivectors whose blades intersect in a 1-dimensional
subspace of $T_{m'}M$ and which again results in condition $(i)$
failing. If dim($\ker\tilde{f})=3$ it can be checked that either
$\ker\tilde{f}$ is spanned by three simple bivectors whose blades
intersect in a 1-dimensional subspace of $T_{m'}M$ (in which case
$d\sigma$ is forced to be zero) or $\ker\tilde{f}$ is spanned by
three simple bivectors the blades of whose duals intersect in a
1-dimensional subspace of $T_{m'}M$ (and $(i)$ fails). If dim
$\ker\tilde{f}\geq 4$, $\ker\tilde{f}$ must contain a non-simple
member and $d\sigma=0$ must hold (see \cite[p.~392]{2}).
 It follows that $d\sigma \equiv0$ on $U$ and, since $U$
is connected, the result $(ii)(a)$ follows. The result $(ii)(b)$ is
then immediate from~(\ref{E25}) and the result $(ii)(c)$ follows
from (\ref{E22}).
\end{proof}

It is remarked that, in the construction of examples, the concept of
{\em local  projective relatedness} will be required. For a
space-time $(M,g)$ let $U$ be a non-empty connected open subset of
$M$ and let $g'$ be some metric def\/ined on $U$. Then $g$ and $g'$
(or their respective Levi-Civita connections) will be said to be
{\em $($locally$)$ projectively related} (on $U$) if the restriction of
$g$ to $U$ is projectively related to~$g'$ on~$U$.

\section{Projective structure and holonomy I}\label{section6}

The relationship between holonomy type and projective relatedness
can now be given. Amongst the holonomy types studied are several
which include signif\/icant solutions to Einstein's f\/ield equations in
general relativity (in addition to the vacuum solutions already
discussed). For example, one has the non-vacuum pp-waves (type
$R_{3}$ and $R_{8}$), the G\"{o}del metric ($R_{10}$), the
Bertotti--Robinson metrics ($R_{7}$) and the Einstein static universe
($R_{13}$). For each holonomy type the general idea is, f\/irst, to
determine the holonomy invariant distributions peculiar to that type
and identify any covariantly constant or recurrent vector f\/ields,
second to link these vector f\/ields to $\rm Riem$ using the inf\/initesimal
holonomy structure, third to use Lemma~\ref{lemma1} and Theorem~\ref{theorem1} to f\/ind
expressions for the 1-form $\lambda$ and the tensor $a$ and f\/inally
to use (\ref{E18}) to complete the procedure. It is convenient to
break up the holonomy types into certain subcollections for easier
handling. All metric and connection statements and index raising are
understood to apply to the structures $g$ and $\nabla$ originally
given on~$M$. The following preliminary topological lemma, which is
a generalised version of a result in~\cite{2}, is useful in some of
the theorems.

\begin{lemma}\label{lemma2}\qquad{}
\begin{enumerate}[$(i)$]\itemsep=0pt
\item Let $X$ be a topological space and let $A$ and $B$ be disjoint subsets of $X$ such that $A$ and $A\cup B$
are open in $X$ and $A\cup B$ is dense in $X$. Suppose $B=B_{1}\cup
B_{2}$, with $B_{1}$ and $B_{2}$ disjoint, and $\mathrm{int}\, B\subset
\mathrm{int}\, B_{1}\cup B_{2}$ where $\mathrm{int}$ denotes the
interior operator in the manifold topology of M. Then $X$ may be
\emph{disjointly} decomposed as $X=A\cup \mathrm{int}\, B_{1}\cup
\mathrm{int}\, B_{2}\cup J$ where $J$ is the closed subset of $X$
defined by the disjointness of the decomposition and $A\cup
\mathrm{int}\,B_{1}\cup \mathrm{int}\,B_{2}$ is open and dense in $X$
$($that is, $\mathrm{int} J={\varnothing})$.

\item Let $X$ be topological space and let $A_{1}, \dots,A_{n}$ be disjoint subsets of $X$ such that $A_{1}$
together with $\cup_{i=1}^{i=k}{A_{i}}$ for $k=2,\dots,n$ are open
subsets of $X$ and such that $\cup_{i=1}^{i=n}{A_{i}}$ is $($open and$)$
dense in $X$. Then $X$ may be \emph{disjointly} decomposed as
$X=A_{1}\cup \mathrm{int}\, A_{2}\cup\cdots\cup \mathrm{int}\,A_{n}\cup K$
where $K$ is the closed subset of $X$ defined by the disjointness of
the decomposition and $A_{1}\cup \mathrm{int}\, A_{2}\cup\cdots \cup
\mathrm{int}\,A_{n}$ is open and dense in $X$ $($that is,
$\mathrm{int}\,K={\varnothing})$.
\end{enumerate}
\end{lemma}

\begin{proof}
$(i)$
Suppose $\mathrm{int}\, J\neq {\varnothing}$ and let $U$ be a non-empty open subset of $X$ with $U\subset J$. Then $U\cap (A\cup B)$
is open and non-empty (since $A\cup B$ is dense in $M$) but $U$ is
disjoint from $A$, $\mathrm{int}\,B_{1}$ and $\mathrm{int}\,B_{2}$. Then
$U\cap B$ is open and non-empty and hence so is $U\cap \mathrm{int}\,
B$. But ${\varnothing} \neq U\cap \mathrm{int}\,B \subset U\cap
(\mathrm{int}\,B_{1}\cup B_{2}) = U\cap B_{2}$ and so $U\cap B_{2}$
contains the non-empty, open subset $U\cap \mathrm{int}\,B$ of $M$
from which the contradiction that $U\cap \mathrm{int}\,B_{2}\neq
{\varnothing}$ follows. Thus, $U={\varnothing}$,
$\mathrm{int}J={\varnothing}$ and the result follows.

$(ii)$ Suppose $U$ is an non-empty open subset of $X$ with $U\subset K$.
Then $U\cap (\cup_{i=1}^{i=n}{A_{i}})$ is not empty but $U$ is
disjoint from $A_{1}$, $\mathrm{int}\,A_{2},\dots,\mathrm{int}\,A_{n}$.
It follows that if the open set $U\cap (A_{1}\cup A_{2}) (=U\cap
A_{2}$) is non-empty, then $U\cap A_{2}$ is non-empty and open and
gives the contradiction that $U\cap \mathrm{int}\,A_{2}$ is (open and)
non-empty. Thus $U$ is disjoint from $A_{2}$. Continuing this
sequence one f\/inally gets the contradiction that $U$ is disjoint
from each $A_{i}$ and hence from the open dense set
$\cup_{i=1}^{i=n}{A_{i}}$. Thus $U={\varnothing}$ and the result
follows.
\end{proof}

\begin{theorem}\label{theorem2}
Let $(M,g)$ and $(M,g')$ be space-times with $(M,g)$ non-flat.
Suppose that $(M,g)$ is of holonomy type  $R_{2}$, $R_{3}$ or
$R_{4}$ and that $\nabla$ and $\nabla'$ are projectively related.
Then $\nabla=\nabla'$ on~$M$.
\end{theorem}
\begin{proof} Suppose f\/irst that $(M,g)$ has holonomy type $R_{2}$ and, by the non-f\/lat condition,
let $U$ be the open dense subset of $M$ on which $\rm Riem$ is non-zero.
Then, for $m \in U$ (see Section~\ref{section4}), there exists a connected and
simply connected open neighbourhood $V\subset U$ of $m$ and two
orthogonal smooth unit spacelike vector f\/ields $X$ and $Y$ on $V$,
spanning a holonomy invariant distribution at each point of $V$ and
which are covariantly constant on $V$. The Ricci identity then
reveals that $R^{a}{}_{bcd}X^{d}=R^{a}{}_{bcd}Y^{d}=0$ on $V$ and so
$\rm Riem$ takes the curvature class $\mathbf{D}$ form~(\ref{E5}) where
$F$ is a smooth simple timelike bivector f\/ield on $V$ whose blade is
orthogonal to the 2-spaces $X(m)\wedge Y(m)$ at each $m\in V$. Then
at each $m\in V$ one may construct a null tetrad $l$, $n$, $x$, $y$, based on
the holonomy invariant subspaces, so that $x=X(m)$ and $y=Y(m)$ and
then $R_{abcd}G^{cd}=0$ is satisf\/ied by at $m$ for $G=l\wedge x$,
$l\wedge y$, $n\wedge x$, $n\wedge y$ and $x\wedge y$. It follows
that the conditions of Lemma~\ref{lemma1} are satisf\/ied on (some possibly
reduced version of) $V$ and hence that the conclusions $(a)$, $(b)$ and
$(c)$ of this lemma hold. Then part $(c)$ of Lemma~\ref{lemma1} and Theorem~\ref{theorem1}$(i)$
show that, on $V$,
\begin{equation}\label{27}
a_{ab}=\phi g_{ab}+\mu X_{a}X_{b}+\nu Y_{a}Y_{b}+\rho (X_{a}Y_{b}+Y_{a}X_{b})
\end{equation}
for functions $\phi$, $\mu$, $\nu$ and $\rho$. Also the smoothness of the functions $a_{ab}X^{a}X^{b}$ $(=\phi+\mu)$,
$a_{ab}Y^{a}Y^{b}$ $(=\phi+\nu)$, $a_{ab}X^{a}Y^{b}$ $(=\rho)$ and $g^{ab}a_{ab}$ $(=4\phi+\mu+\nu)$ reveal the smoothness of the functions
$\phi$, $\mu$, $\nu$ and $ \rho$ on $V$. Now one substitutes (\ref{27}) into (\ref{E18}) to get
\begin{equation}\label{28} g_{ab}\phi_{,c}+X_{a}X_{b}\mu_{,c}+Y_{a}Y_{b}\nu_{,c}+(X_{a}Y_{b}+Y_{a}X_{b})\rho_{,c}=g_{ac}\lambda_{b}+ g_{bc}\lambda_{a}.
 \end{equation}
A contraction of (\ref{28}) successively with $l^{a}x^{b}$ and
$n^{a}y^{b}$ at any $m\in V$ shows that
$\lambda_{a}x^{a}=\lambda_{a}l^{a}=\lambda_{a}y^{a}=\lambda_{a}n^{a}=0$.
Thus the 1-form $\lambda$ is zero on $V$. It follows from
(\ref{E17}) that the 1-form $\psi$ is identically zero on $V$ and
hence, from (\ref{E12}) that $\nabla=\nabla'$ on $V$, hence on $U$
and thus on~$M$. This completes the proof for holonomy type $R_{2}$
and the proofs for the holonomy types $R_{3}$ and $R_{4}$ are
similar; for type $R_{3}$ one has covariantly constant, orthogonal,
null and spacelike vector f\/ields $L$ and $Y$, respectively, on $V$
and the holonomy invariant distributions can be used construct
vector f\/ields $N$ and $X$ on (a possibly reduced) $V$ such that $L,
N, X$ and $Y$ give a null tetrad at each point of $V$ from which the
proof follows in a similar way to that of the $R_{2}$ type. For type
$R_{4}$ one has covariantly constant null vector f\/ields $L$ and $N$
on $V$ which one can choose to satisfy $g_{ab}L^{a}N^{b}=1$. Again
one easily achieves $\nabla=\nabla'$ on $M$.
\end{proof}

\begin{theorem}\label{theorem3}
Let $(M,g)$ and $(M,g')$ be space-times with $(M,g)$ non-flat and of
holonomy type~$R_{7}$. Suppose that $\nabla$ and $\nabla'$  are
projectively related. Then $\nabla=\nabla'$ on~$M$.
\end{theorem}
\begin{proof}
Let $U$ be the open dense subset of $M$ on which $\rm Riem$ does not
vanish. If $m\in U$ there exists an open, connected and simply
connected neighbourhood $V\subset U$ and null properly recurrent
vector f\/ields $L$ and $N$ scaled so that $L^{a}N_{a}=1$ on $V$
(noting that this scaling will not af\/fect their recurrence
property). Then $\nabla L=L\otimes P$ and $\nabla N=-N\otimes P$ for
some smooth 1-form f\/ield $P$ on $V$. The simple bivector f\/ield
$F\equiv L\wedge N$, ($F_{ab}\equiv 2L_{[a}N_{b]}$), then satisf\/ies
$\nabla F=0$ and its (simple) dual bivector f\/ield also satisf\/ies
$\nabla \overset{*}{F}=0$. The blades of $F$ and $\overset{*}{F}$
span, respectively, the timelike and spacelike holonomy
distributions on $V$ and reducing $V$, if necessary, one may write
$\overset{*}{F}_{ab}=2X_{a[}Y_{b]}$ for smooth unit orthogonal
spacelike vector f\/ields $X$ and $Y$ on $V$ and which satisfy
$X^{a}X_{a;b}=Y^{a}Y_{a;b}=0$. The condition $\nabla
\overset{*}{F}=0$ on $V$ shows that $\nabla X=Y\otimes Q$ and
$\nabla Y=-X\otimes Q$ on $V$ for some smooth 1-form $Q$ on $V$. At
each $m\in V$, $L(m)$, $N(m)$, $X(m)$ and $Y(m)$ give a null tetrad
$l$, $n$, $x$, $y$ at~$m$. Now at each $m\in U$, Table~\ref{table1} and a consideration
of the inf\/initesimal holonomy structure of $(M,g)$ show that the
curvature class of $\rm Riem$ at $m$ is either $\mathbf{D}$ or
$\mathbf{B}$, taking the form (\ref{E4}) if it is $\mathbf{B}$ and
(\ref{E5}) with $F$ either timelike or spacelike if it is
$\mathbf{D}$. Let $\mathbf{B}$, $\mathbf{D_{s}}$ and
$\mathbf{D_{t}}$ denote the subsets of $U$ consisting of those
points at which the curvature class is, respectively, $\mathbf{B}$,
$\mathbf{D}$ with $F$ spacelike or $\mathbf{D}$ with $F$ timelike.
Then one may decompose~$M$, disjointly, as
$M=\mathbf{B}\cup\mathbf{D_{s}} \cup \mathbf{D_{t}}\cup
J=\mathbf{B}\cup \mathrm{int}\,\mathbf{D_{s}} \cup \mathrm{int}\,
\mathbf{D_{t}}\cup K$ where, by the ``rank'' theorem (see e.g.~\cite{2}), $\mathbf{B}$ is open in $M$ and $J$ and $K$ are closed subsets
of $M$ def\/ined by the disjointness of the decomposition and $J$
satisf\/ies $\mathrm{int}\,J$=${\varnothing}$. If
$\mathrm{int}\mathbf{D_{s}}\neq {\varnothing}$ let $m\in \mathrm{int}\,
\mathbf{D_{s}}$ and $m\in V \subset \mathrm{int}\,\mathbf{D_{s}}$
 with $V$ as above. Then (\ref{E5}) holds on $V$ with $F=X\wedge Y$ and $\alpha$ smooth on $V$. It is then clear that,
at $m$, $l\wedge x$, $l\wedge y$, $n\wedge x$, $n\wedge y$ and $l\wedge n$
are in $\ker\tilde{f}$ and so, from Lemma~\ref{lemma1}$(ii)(a)$, $T_mM$ is an eigenspace of $\nabla\lambda$ and so $\nabla
\lambda=cg$ on $V$. Then Lemma~\ref{lemma1}$(ii)(c)$ and Theorem~\ref{theorem1}$(i)$ show
that
\begin{equation}\label {29} a_{ab}=\phi g_{ab}+\mu L_{a}L_{b}+\nu
N_{a}N_{b}+\rho (L_{a}N_{b}+N_{a}L_{b}) \end{equation} for functions
$\phi$, $\mu$, $\nu$, and $\rho$ on $V$ which are smooth since
$a_{ab}L^{a}L^{b}$, $a_{ab}N^{a}N^{b}$, $a_{ab}L^{a}N^{b}$ and
$g^{ab}a_{ab}$ are. Noting that $L_{(a}N_{b)}$ is covariantly
constant on~$V$, a substitution of (\ref{29}) into (\ref{E18})
gives, on~$V$,
 \begin{gather}\label{30}
 g_{ab}\phi_{,c}+L_{a}L_{b}\mu_{,c}+2\mu L_{a}L_{b}P_{c}+N_{a}N_{b}\nu _{,c}-2\nu N_{a}N_{b}P_{c}+(L_{a}N_{b}+N_{a}L_{b})\rho _{,c}\\\nonumber
 \qquad {}=g_{ac}\lambda _{b}+g_{bc}\lambda_{a}.
 \end{gather}
Successive contractions of (\ref{30}) with $L^{a}X^{b}$ and
$N^{a}Y^{b}$ then show that $\lambda_{a}L^{a}=
\lambda_{a}N^{a}=\lambda_{a}X^{a}=\lambda_{a}Y^{a}=0$ on $V$. Thus
$\lambda$ vanishes on $V$ and so on $\mathrm{int}\,\mathbf{D_{s}}$. If
$\mathrm{int}\,\mathbf{D_{t}}\neq {\varnothing}$ a similar argument
using~(\ref{E5}) with $F=L\wedge N$ can be used to show that
$\lambda$ vanishes on $\mathrm{int}\,\mathbf{D_{t}}$. If
$\mathbf{B}\neq {\varnothing}$, Lemma~\ref{lemma1}$(ii)(b)$ shows that $\lambda$
vanishes on $\mathbf{B}$ since no non-trivial solutions for
$\lambda$ of the equation displayed there exist at points of
curvature class $\mathbf{B}$. So one has achieved the situation that
$M=\mathbf{B}\cup \mathbf{D}\cup J$ with $\mathrm{int}\,J={\varnothing}$
and $\mathbf{B}$ and $\mathbf{B}\cup \mathbf{D}$ open in $M$ (again
by an application of the ``rank'' theorem) and with $\lambda$
vanishing on the open subset $\mathbf{B}\cup
\mathrm{int}\,\mathbf{D_{s}}\cup \mathrm{int}\,\mathbf{D_{t}}$. That
this open subset is a dense subset of $M$ follows from Lemma~\ref{lemma2}$(i)$ by
the following argument. Since $\mathbf{D}\equiv \mathbf{D_{s}}\cup
\mathbf{D_{t}}$ then, from elementary topology,
$\mathrm{int}\,\mathbf{D_{s}}\cup \mathrm{int}\,\mathbf{D_{t}}\subset
\mathrm{int}\,\mathbf{D}$. That the reverse inclusion is true in this
case follows from noting that it is trivially true if
$\mathrm{int}\,\mathbf{D}={\varnothing}$ whilst if
$\mathrm{int}\,\mathbf{D}\neq {\varnothing}$ one writes equation
(\ref{E5}) for $\rm Riem$ on some open neighbourhood of any point $m$ in
the submanifold $\mathrm{int}\,\mathbf{D}$ with $\alpha$ and $F$
smooth (as one can) and then def\/ines the smooth map $\theta$ from
this neighbourhood to $\mathbb{R}$ by $\theta (m)=(F^{ab}F_{ab})(m)$
(so that $\theta (m)> 0$ if $m\in \mathbf{D_{s}}$ and $\theta (m)<
0$ if $m\in \mathbf{D_{t}}$). It is easily seen now that $m$ is
either in $\mathrm{int}\,\mathbf{D_{s}}$ or
$\mathrm{int}\,\mathbf{D_{t}}$. Thus
$\mathrm{int}\,\mathbf{D}=\mathrm{int}\,\mathbf{D_{s}}\cup
\mathrm{int}\,\mathbf{D_{t}}$ and the conditions of Lemma~\ref{lemma2}$(i)$ are
satisf\/ied. So $\lambda$ vanishes on an open dense subset of $M$ and
hence on~$M$. It follows that $\nabla=\nabla'$ on~$M$.
\end{proof}

\begin{theorem}\label{theorem4}
Let $(M,g)$ and $(M,g')$ be space-times with $(M,g)$ non-flat and of
holonomy  type~$R_{10}$, $R_{11}$ or $R_{13}$ and with $\rm Riem$ of
curvature class $\mathbf{C}$ at each point of an open, dense subset~$U$ of~$M$. Suppose that $\nabla$ and $\nabla'$ are projectively
related. Then $\nabla=\nabla'$ on~$M$.
\end{theorem}

\begin{proof}
Suppose that $(M,g)$ is of holonomy type $R_{10}$ (the arguments for
types $R_{11}$ and $R_{13}$ are similar). For any $m\in U$ there
exists an open, connected and simply connected neighbourhood
$V\subset U$ of $m$ and a unit spacelike covariantly constant vector
f\/ield $X$ on $V$. Then one may choose an orthonormal tetrad
$u$, $x$, $y$, $z$ at $m$ such that $X(m)=x$ and, from the Ricci identity on
$X$, $R^{a}{}_{bcd}X^{d}=0$ on $V$, and so $x\wedge u$, $x\wedge y$
and $x\wedge z$ are in $\ker\tilde{f}$ at $m$. As before, one sees
that the conditions and conclusions of Lemma~\ref{lemma1} hold at each point of
$V$ and so Theorem~\ref{theorem1}$(ii)$ gives
\begin{equation}\label{31}
a_{ab}=\phi g_{ab}+\mu X_{a}X_{b}
\end{equation}
on $V$ for functions $\phi$ and $\mu$ which are smooth since
$a_{ab}X^{a}X^{b}$ and $g^{ab}a_{ab}$ are. On substituting~(\ref{31}) into (\ref{E18}) one f\/inds
\begin{equation}\label{32}
g_{ab}\phi_{,c}+X_{a}X_{b}\mu _{,c}=g_{ac}\lambda_{b}+g_{bc}\lambda_{a}
\end{equation}
and successive contractions of (\ref{32}) with $x^{a}y^{b}$ and
$u^{a}z^{b}$ show that $\lambda$ vanishes at $m$, hence on $V$ and
so on $M$. It follows that $\nabla=\nabla'$ on $M$.
\end{proof}

\begin{theorem}\label{theorem5}
Let $(M,g)$ and $(M,g')$ be space-times with $(M,g)$ non-flat and of
holonomy type~$R_{6}$, $R_{8}$ or $R_{12}$.  Suppose that $\nabla$ and
$\nabla'$ are projectively related. Then $\nabla=\nabla'$ on $M$.
\end{theorem}

\begin{proof}
In each case let $U$ be the open dense subset of $M$ on which $\rm Riem$
does not vanish. Then any $m\in U$ admits an open, connected and
simply connected neighbourhood $V\subset U$ on which are def\/ined
orthogonal vector f\/ields $L$, $X$ and $Y$, with $L$ null and either
recurrent or covariantly constant (depending on the holonomy type)
and $X$ and $Y$ unit spacelike vector f\/ields and which together span
a 3-dimensional null holonomy invariant distribution on $V$. By
reducing $V$, if necessary, one may assume the existence of a null
vector f\/ield $N$ which along with $L$, $X$ and $Y$ span a null
tetrad at each point of $V$.

Now for holonomy type $R_{6}$ one may take $L$ recurrent and $Y$
covariantly constant. The curvature class at $m$ is either
$\mathbf{C}$ or $\mathbf{D}$ and if it is $\mathbf{D}$, equation
(\ref{E5}) holds with $F=l\wedge x$ or $F=l\wedge n +al\wedge x$
($a\in\mathbb{R}$). Let the subset of points of $M$ at which the
f\/irst of these class $\mathbf{D}$ possibilities holds be denoted by
$\mathbf{D_{n}}$ (null) and at which the second holds by
$\mathbf{D_{nn}}$ (non-null). Then let $\mathbf{D}\equiv
\mathbf{D_{n}}\cup \mathbf{D_{nn}}$. Let the subset of points of $M$
at which at which the curvature class is~$\mathbf{C}$ be denoted by
$\mathbf{C}$. Then $\mathbf{C}$ is open in $M$ and $U=\mathbf{C}\cup
\mathbf{D}$. $M$ admits the disjoint decomposition $M=\mathbf{C}\cup
\mathbf{D_{n}} \cup \mathbf{D_{nn}} \cup J$ where $J$ is the closed
subset of $M$ determined by the disjointness of the decomposition
and $\mathrm{int}\, J={\varnothing}$. It is clear that, whatever the
curvature type, the conditions and conclusions of Lemma~\ref{lemma1} hold. If
$\mathbf{C}\neq \varnothing$ and $m\in \mathbf{C}$ then one may
arrange in the previous paragraph that $m\in V \subset U$ and that,
on $V$,
\begin{equation}\label{33}
a_{ab}=\phi g_{ab}+\mu Y_{a}Y_{b}
\end{equation}
for functions $\phi$ and $\mu$. The proof now proceeds as for
Theorem~\ref{theorem4}  and one obtains $\nabla=\nabla'$ on~$\mathbf{C}$.

If $\mathrm{int}\,\mathbf{D_{n}}\neq \varnothing$ then with $V$ as
above and $m\in V\subset \mathrm{int}\,\mathbf{D_{n}}$, Lemma~\ref{lemma1} and
Theorem~\ref{theorem1}$(i)$ give, on~$V$,
\begin{equation}\label{34}
a_{ab}=\phi g_{ab}+\mu L_{a}L_{b}+\nu Y_{a}Y_{b}+\rho (L_{a}Y_{b}+Y_{a}L_{b})
\end{equation} for functions $\phi$, $\mu$, $\nu$ and $\rho$ which are smooth since $a_{ab}X^{a}X^{b}$, $a_{ab}Y^{a}Y^{b}$, $a_{ab}N^{a}N^{b}$ and
$a_{ab}N^{a}Y^{b}$ are. On substituting (\ref{34}) into (\ref{E18})
and contracting, successively, with $X^{a}Y^{b}$, $L^{a}L^{b}$ and
$N^{a}X^{b}$ one again f\/inds that $\lambda=0$ on $V$ and hence on
$\mathrm{int}\mathbf{D_{n}}$ and so $\nabla=\nabla'$ on
$\mathrm{int}\mathbf{D_{n}}$. (This argument includes within it the
curvature class $\mathbf{C}$ case above; one simply sets $\mu$ and
$\rho$ to zero in~(\ref{34}).)

If  $\mathrm{int}\,\mathbf{D_{nn}}\neq \varnothing$ and $m\in V\subset
\mathrm{int}\mathbf{D_{nn}}$ with $V$ as above, $\rm Riem$ takes the
form (\ref{E5}) with $F$ chosen as $F=L\wedge N+\sigma L\wedge X$
for a smooth function $\sigma$, on $V$. Now def\/ine smooth vector
f\/ields $X'=X-\sigma L$ and $N'=N+\sigma X-\frac{\sigma ^{2}}{2}L$ on
$V$, noting that $L$, $N'$, $X'$ and $Y$ still give rise to a null
tetrad in the obvious way at each point of $V$ and that now
$F=L\wedge N'$ on $V$. Again Lemma~\ref{lemma1} applies and Theorem~\ref{theorem1}$(i)$ gives
\begin{equation}\label{35} a_{ab}=\phi g_{ab}+\mu X'_{a}X'_{b}+\nu
Y_{a}Y_{b}+\rho (X'_{a}Y_{b}+Y_{a}X'_{b}),
\end{equation}
where, as before, $\phi$, $\mu$, $\nu$ and $\rho$ are smooth functions
on $V$.  On substituting (\ref{35}) into (\ref{E18}) (and noting
that, since $L$ is recurrent, $L^{a}X'_{a;b}=0)$ successive
contractions with $N'^{a}N'^{b}$, $L^{a}X'^{b}$ and $L^{a}Y^{b}$ again
reveal that $\nabla=\nabla'$ on $\mathrm{int}\,\mathbf{D_{nn}}$. So
$\nabla=\nabla'$ on $\tilde{W}\equiv \mathbf{C}\cup
\mathrm{int}\,\mathbf{D_{n}} \cup \mathrm{int}\,\mathbf{D_{nn}}$. That
the open set $\tilde{W}$ is dense in $M$ follows from Lemma~\ref{lemma2}$(i)$. To
see this let $K\equiv M\setminus \tilde{W}$ and let~$W$ be a
non-empty open subset of $K$. Then $W$ is disjoint from
$\mathbf{C}$, $\mathrm{int}\,\mathbf{D_{n}}$ and
$\mathrm{int}\,\mathbf{D_{nn}}$ but the open subset $W\cap
(\mathbf{C}\cup \mathbf{D})=W\cap U \neq {\varnothing}$ since $U$ is
dense in $M$. It follows that $W\cap \mathbf{D}$ and hence $W\cap
\mathrm{int}\mathbf{D}$ are open and non-empty. Now, as in the
previous theorem, consider the map arising from the bivector $F$ in
(\ref{E5}) on the (non-empty) open submanifold
$\mathrm{int}\mathbf{D}$ (and note that for points in
$\mathbf{D_{n}}$ (where $F$ null) its value is zero whereas for
points in $\mathbf{D_{nn}}$ its value is negative). A consideration
of the open subset of $\mathrm{int}\mathbf{D}$ given by the inverse
image, under $f$, of the negative reals shows that
$\mathrm{int}\mathbf{D}\subset \mathrm{int}\mathbf{D_{nn}}\cup
\mathbf{D_{n}}$. The result now follows from Lemma~\ref{lemma2}$(i)$ and
$\lambda$ vanishes on $\tilde{W}$ and hence on~$M$. It follows that
$\nabla=\nabla'$ and this completes the proof when the holonomy type
is $R_{6}$.

If $(M,g)$ has holonomy type $R_{8}$ and if $m\in U$ the curvature
class of $\rm Riem$ at $m$ is either $\mathbf{D}$ (equation (\ref{E5})
with $F$ null) or $\mathbf{C}$. A neighbourhood  $V$ and vector
f\/ields $L$, $N$, $X$ and $Y$, can be established, as in the above
paragraphs, with $L$ covariantly constant.  Lemma 1 again applies
and Theorem~\ref{theorem1}, parts $(i)$ and $(ii)$, then lead to a straightforward
proof that $\lambda=0$ on (and in an obvious notation) the open
subset $\mathbf{C}$ and on $\mathrm{int}\,\mathbf{D}$. (For class
$\mathbf{D}$, useful contractions of (\ref{E18}) are with
$L^{a}L^{b}$, $L^{a}X^{b}$, $L^{a}Y^{b}$, $Y^{a}Y^{b}$ and $L^{a}N^{b}$
and in the class $\mathbf{C}$ case, with $L^{a}X^{b}$, $L^{a}Y^{b}$
and $X^{a}N^{b}$.) That the open subset $\mathbf{C}\cup
\mathrm{int}\,\mathbf{D}$ is dense in $M$ then follows from Lemma~\ref{lemma2}$(ii)$ since $\mathbf{C}\cup \mathbf{D}$ is open and dense in $M$.
Thus $\lambda=0$ and $\nabla=\nabla'$ on $M$.

If $(M,g)$ has holonomy type $R_{12}$ then for $m\in U$ one may
choose a neighbourhood $V$ of $m$ and vector f\/ields $L$, $N$, $X$ and
$Y$ as above with the null vector f\/ield $L$ recurrent. Table~\ref{table1} shows
that $l\wedge x$, $l\wedge y$ and $l\wedge n+\omega ^{-1}x\wedge y$
lie in $\ker\tilde{f}$ and so Lemma~\ref{lemma1} again applies. If $m\in U$ the
possible curvature classes for $\rm Riem$ at $m$ are $\mathbf{D}$ with
$F$ null in~(\ref{E5}) (since the only linear combinations of the
members of the bivector algebra for type $R_{12}$ which are simple
are null), $\mathbf{C}$ (with $R^{a}{}_{bcd}l^{d}=0$ at $m$) or
$\mathbf{A}$. One then decomposes~$M'$, in an obvious notation, as
$M'=\mathbf{A}\cup \mathrm{int}\,\mathbf{C}\cup
\mathrm{int}\,\mathbf{D}\cup J$ (recalling from Section~\ref{section3} that
$\mathbf{A}$ is an open subset of~$M$). The conclusion $(b)$ of Lemma~\ref{lemma1}$(ii)$ then immediately shows that $\lambda=0$ on $\mathbf{A}$ and
procedures similar to those already given reveal that $\lambda=0$ on
$\mathrm{int}\,\mathbf{C}$ and $\mathrm{int}\,\mathbf{D}$ (for the
latter contract~(\ref{E18}) with $L^{a}L^{b}$, $L^{a}X^{b}$,
$L^{a}Y^{b}$, $Y^{a}Y^{b}$ and $N^{a}L^{b}$). Since $\mathbf{A}$,
$\mathbf{A}\cup \mathbf{C}$ and $\mathbf{A}\cup \mathbf{C}\cup
\mathbf{D}$ are open in $M$, with the latter dense, Lemma~\ref{lemma2}$(ii)$
completes the proof that $\lambda=0$ and $\nabla=\nabla'$ on~$M$.
\end{proof}

In summary, if $(M,g)$ is non-f\/lat and if, in addition, it is either
of  holonomy type $R_{2}$, $R_{3}$, $R_{4}$, $R_{6}$, $R_{7}$,
$R_{8}$ or $R_{12}$ or of holonomy type $R_{10}$, $R_{11}$ or $R_{13}$
\emph{and} of curvature class $\mathbf{C}$ over some open dense
subset of~$M$ then if $g'$ is another metric on~$M$ projectively
related to $g$, $\nabla=\nabla'$ on~$M$. One can say a little more
here because in each of the above holonomy types (including the
curvature class $\mathbf{C}$ clause applied to the types
$R_{10}$, $R_{11}$ and $R_{13}$) the condition that $\nabla$ and
$\nabla'$ are projectively related leads to $\nabla=\nabla'$ and
hence to the condition that the holonomy groups of $(M,g)$ and
$(M,g')$ are the same. It is then easy to write down a simple
relation between the metrics $g$ and $g'$~\cite{14, 2}. For example,
suppose either that $M$ is simply connected or that one is working
over a connected and simply connected open region $V$ of $M$ (on
which the holonomy type is the same as that of $M$). Then if this
common holonomy type is $R_{2}$, $R_{3}$ or $R_{4}$, one gets a
relation of the form
\begin{equation}\label{36}
 g'_{ab}=\phi g_{ab}+\mu p_{a}p_{b}+\nu q_{a}q_{b}+\rho (p_{a}q_{b}+q_{a}p_{b}),
 \end{equation}
 where $\phi$, $\mu$, $\nu$ and $\rho$ are constants (subject only to a non-degeneracy condition and signature requirements) and where $p$ and $q$ are vector f\/ields on $M$ (or $V$) spanning the 2-dimensional vector space of covariantly constant vector f\/ields. For holonomy types $R_{6}$ and $R_{8}$ or types $R_{10}$, $R_{11}$ or $R_{13}$ (with the above clauses attached) one gets
 \begin{equation}\label{37}
 g'_{ab}=\phi g_{ab}+\mu k_{a}k_{b},
 \end{equation}
 where $\phi$ and $\mu$ are constants (again subject to non-degeneracy and signature requirements) and~$k$ is a vector f\/ield on $M$ (or $V$) spanning the vector space of covariantly constant vector f\/ields. For holonomy type $R_{7}$ one gets
\begin{equation}\label{38}
g'_{ab}=\phi g_{ab}+\mu (l_{a}n_{b}+n_{a}l_{b}),
\end{equation}
where again $\phi$ and $\mu$ are constants (again subject to
non-degeneracy and signature requirements) and $l$ and $n$ are null,
properly recurrent vector f\/ields on $M$ (or $V$). For holonomy type
$R_{12}$ one gets $g'_{ab}=\phi g_{ab}$ with $\phi$ constant. [It is
remarked here that for holonomy types $R_{10}$, $R_{11}$ or
$R_{13}$, \emph{even without the curvature class $\mathbf{C}$
clause, if it is given that $\nabla=\nabla'$ then \eqref{37} follows
for constant $\phi$ and $\mu$}. But, in this latter case, projective
relatedness without this clause does not force the condition
$\nabla=\nabla'$ as will be seen later.] It is clear from these
results that if $g$ and $g'$ are projectively related metrics they
may have dif\/ferent signatures.

\section{Projective structure and holonomy II}\label{section7}

In this section the holonomy types (and special cases of holonomy
types) omitted in the last section will be discussed.  Theorem~\ref{theorem4}
dealt in detail with the projective relatedness problem for the
(non-f\/lat) holonomy types $R_{10}$, $R_{11}$ and $R_{13}$ when the
curvature rank of the map $f$ in~(\ref{E3}) (and which, for each of
these types, is necessarily at most~3) is equal to~2 or~3 (curvature
class~$\mathbf{C}$) at each point of an open, dense subset of~$M$.
One can now consider the situation for these holonomy types when the
curvature rank is~1 (curvature class $\mathbf{D}$) at each point of
some open, dense subset $U$ of $M$. For any of these holonomy types
def\/ine subsets $M_{1} \equiv \{m\in M : \lambda$ vanishes on some
open neighbourhood of $m\}$ and $M_{2} \equiv \{m\in M: \lambda$
does not vanish on any open neighbourhood of $m\} = M \setminus
M_{1}$. Then one has the disjoint decomposition of $M$ given by $M =
M_{1}\cup M_{2}$ with $M_{1}$ open in $M$. Let $J=M\setminus U$ so
that $J$ is closed in $M$ and $\mathrm{int}\,J=\varnothing$ (since $U$ is open
and dense in $M$). Next, if $\mathbf{D_{nn}}$ is the subset of points of M at
which the bivector~$F$ in~(\ref{E5}) is either spacelike or timelike and
$\mathbf{D_{n}}$ the subset of points where it is null then $U=\mathbf{D_{nn}}\cup
\mathbf{D_{n}}$ with $\mathbf{D_{nn}}$ an open subset of $M$ (by a consideration of
the map associated with the smooth function $F_{ab}F^{ab}$ on $U$
and described earlier). One then gets a disjoint decomposition $M =
M_{1} \cup (M_{2} \cap \mathbf{D_{nn}}) \cup\mathrm{int}\,(M_{2} \cap \mathbf{D_{n}}) \cup K$.
It will be seen, as the proof progresses, that in all cases $M_{2}
\cap \mathbf{D_{nn}}$ is also open in $M$ and hence $K$, which is def\/ined by
the disjointness of the decomposition, is closed in $M$. It then
follows that, with the openness of $M_{2} \cap \mathbf{D_{nn}}$ assumed, if~$W$ is a non-empty open subset of $K$, $W \cap M_{1} = W \cap (M_{2}
\cap \mathbf{D_{nn}}) =W\cap \mathrm{int}\, (M_{2}\cap \mathbf{D_{n}}) =\varnothing$. But
the non-empty open set $W\cap U=W\cap (\mathbf{D_{n}}\cup \mathbf{D_{nn}})=W\cap
\mathbf{D_{n}}$ (since $W\subset K$ and $K\cap \mathbf{D_{nn}}=\varnothing)$$=W \cap
(M_{2} \cap \mathbf{D_{n}})$. Hence $W\cap \mathrm{int}\, (M_{2} \cap \mathbf{D_{n}})\neq
\varnothing$ and this contradiction shows that $\mathrm{int}\, K = \varnothing$.
Thus, considering the open dense subset $M_{1} \cup (M_{2} \cap
\mathbf{D_{nn}}) \cup\mathrm{int}\,(M_{2}\cap \mathbf{D_{n}})$ of M one sees that since $\lambda = 0$
on $M_{1}$ the interesting parts of $M$ in this context are
the open subsets $M_{2} \cap \mathbf{D_{nn}}$ and $\mathrm{int}\,(M_{2} \cap
\mathbf{D_{n}})$ (and, in fact, $\mathbf{D_{n}} =\varnothing$ in the $R_{13}$ case). In
fact, it will turn out that, in all cases, $\mathrm{int}\,(M_{2}\cap
\mathbf{D_{n}})=\varnothing$.

So let $(M,g)$ be a (non-f\/lat) space-time of holonomy type $R_{10}$,
$R_{11}$ or $R_{13}$ whose curvature is of class $\mathbf{D}$ on
some open dense subset $U$ of $M$. Thus~(\ref{E5}) holds on $U$ (and
$\rm Riem$ vanishes on $M\setminus U$). Let $m\in U$ and let $V\subset
U$ be a connected, simply connected, open neighbourhood of $m$ on
which is def\/ined a metric $g'$, (locally) projectively related to
$g$ on $V$, together with an associated pair $(a, \lambda)$. The
conditions imposed on $(M,g)$ together with (\ref{E5}) show that the
conditions and hence the conclusions of Lemma~\ref{lemma1} hold on $U$ and
hence on $V$. Thus $\nabla \lambda=cg$ holds on $V$ for some
constant $c$.

Now suppose that the holonomy type of $(M,g)$ is $R_{11}$. Then $V$
can be chosen so that it admits a covariantly constant, null vector
f\/ield~$l$, so that, from the Ricci identity, $R^{a}{}_{bcd}l^{d}=0$.
Also, from the curvature class $\mathbf{D}$ condition, $V$ (after a
possible reduction, if necessary) admits another (smooth)
nowhere-zero vector f\/ield $p$ such that $l(m')$ and $p(m')$ are
independent members of $T_{m'}M$ for each $m'\in V$ and such that
$R^{a}{}_{bcd}p^{d}=0$ holds on $V$. Condition $(b)$ of Lemma~\ref{lemma1} and
the curvature class $\mathbf{D}$ condition then show that, on $V$,
\begin{equation}\label{43}
\lambda_{a}=\sigma l_{a}+\rho p_{a}
\end{equation}
and so, since $\lambda_{a;b}=cg_{ab}$, one gets
\begin{equation}\label{44}
\rho p_{a;b}+p_{a}\rho_{,b}+l_{a}\sigma_{,b}=cg_{ab}. \end{equation}
Also, part $(c)$ of Lemma~\ref{lemma1} and Theorem~\ref{theorem1}$(i)$ give
\begin{equation}\label{45}
a_{ab}=\phi g_{ab}+\beta l_{a}l_{b}+\gamma p_{a}p_{b}+\delta (l_{a}p_{b}+p_{a}l_{b})
\end{equation}
for smooth functions $\phi$, $\beta$, $\gamma$ and $\delta$ on $V$. Now
there exists a nowhere-null, nowhere-zero vector f\/ield $q$ on (a
possibly reduced) $V$ which is orthogonal to $l$ and $p$ and is such
that $q_{a}q^{a}$ is constant (and hence $q^{a}q_{a;b}=0$) on $V$.
Then from (\ref{45}),
$a_{ab}q^{b}=\phi q_{a}$ on $V$. On covariantly
dif\/ferentiating this last equation and using (\ref{E18}) after a
contraction with $q^{a}$ one f\/inds that $\phi _{,a}=0$ and so $\phi$
is constant on $V$ (and non-zero since $a$ is non-degenerate).

Now the simple bivector $F$ in (\ref{E5}) satisf\/ies $F_{ab}l^{b}=0$
and so is either null or spacelike at $m\in U$, that is, either
$m\in \mathbf{D_{n}}$ or $m\in \mathbf{D_{nn}}$. Recalling the decomposition of $M$
given above assume $\mathrm{int}(M_{2}\cap \mathbf{D_{n}})$ is not empty and let $m\in
V\subset\mathrm{int}(M_{2}\cap \mathbf{D_{n}})$ for (some possibly reduced) open
neighbourhood $V$ as before. Then $p$ may be chosen a unit spacelike
vector f\/ield on $V$ orthogonal to $l$ on $V$ and (\ref{45}) shows
that $a_{ab}l^{b}=\phi l_{a}$ which, after dif\/ferentiation and
making use of (\ref{E18}) and the constancy of $\phi$ and $l$
immediately shows that $\lambda$ vanishes on $V$. This contradicts
the fact that $m\in M_{2}$. Thus, $\mathrm{int}(M_{2}\cap \mathbf{D_{n}})$ is empty.

If $F$ is spacelike at $m$, that is, $m\in \mathbf{D_{nn}}$, and assuming
that $M_{2}\cap \mathbf{D_{nn}}$ is not empty let $m\in (M_{2}\cap \mathbf{D_{nn}})$
and choose $V$ such that $m\in V\subset \mathbf{D_{nn}}$ (since $\mathbf{D_{nn}}$ is
open in $M$). Then one may choose $p$ so that $p$ is null and
$l^{a}p_{a}=1$ on $V$. Then (\ref{44}) on contraction f\/irst with
$l^{a}$ and then with $p^{a}$ (using $p^{a}p_{a;b}=0$ since $p$ is
null on $V$) yields
\begin{equation}\label{46} cl_{a}=\rho_{,a},\qquad
cp_{a}=\sigma_{,a}
\end{equation} and then (\ref{44}) becomes
\begin{equation}\label{47}
\rho p_{a;b}=cT_{ab},\qquad (T_{ab}\equiv g_{ab}-l_{a}p_{b}-p_{a}l_{b}).
\end{equation}
Next, dif\/ferentiate (\ref{45}) and use (\ref{43}) and (\ref{E18}) to get
\begin{gather}\label{48}
g_{ac}(\sigma l_{b}+\rho p_{b})+g_{bc}(\sigma l_{a}+\rho p_{a})\\\nonumber
\qquad{}=l_{a}l_{b}\beta _{,c}+p_{a}p_{b}\gamma _{,c}+\gamma p_{a;c}p_{b}+\gamma p_{a}p_{b;c}+(l_{a}p_{b}+p_{a}l_{b})\delta _{,c}+\delta (l_{a}p_{b;c}+p_{a;c}l_{b}).
\end{gather}
On contracting this successively with $l^{a}l^{b}$, $p^{a}p^{b}$ and $l^{a}p^{b}$ (using $l^{a}p_{a;b}=0$) one gets
\begin{equation}\label{49}
2\rho l_{a}=\gamma_{,a}, \qquad 2\sigma p_{a}=\beta_{,a},\qquad \sigma l_{a}+\rho p_{a}=\delta _{,a}. \end{equation}
If one f\/inally contracts (\ref{48}) f\/irst with $l^{a}$ and then with $p^{a}$ and use is made of (\ref{49}) one easily f\/inds
\begin{equation}\label{50}
\gamma p_{a;b}=\rho T_{ab}, \qquad \delta p_{a;b}=\sigma T_{ab}.
\end{equation}

The condition $m\in M_{2}$ implies that $p_{a;b}(m)\neq 0$ because,
otherwise, (\ref{47}) shows that $c=0$ and then (\ref{46}) reveals
that $\rho$ and $\sigma$ are constant on $V$. Then (\ref{50}) shows
that $\rho$ and $\sigma$ vanish at~$m$ and hence on $V$. Finally,
(\ref{43}) implies that $\lambda$ vanishes on $V$ and a
contradiction follows. It thus follows that $p_{a;b}$ does not
vanish at~$m$ and hence on some open neighbourhood of $m$ (taken to
be $V$). Further, if $\lambda(m)=0$ then (\ref{43}) shows that
$\sigma(m)=\rho(m)=0$ and so, from (\ref{47}), $c=0$. Then
(\ref{46}) reveals that $\sigma$ and $\rho$ are constant and hence
zero on $V$ and hence, from (\ref{43}), the contradiction that
$\lambda\equiv 0$ on $V$. Thus $\lambda(m)\neq 0$ and so one may
assume that $V$ is chosen (after a possible reduction) so that
$\lambda$ and $p_{a;b}$ are nowhere zero on $V$. But this implies
that $V\subset M_{2}$ and hence that $V\subset (M_{2}\cap \mathbf{D_{nn}})$.
From this it follows that $M_{2}\cap \mathbf{D_{nn}}$ is open in $M$ (as
remarked above).

Continuing with this case, if $p_{a;b}$ is \emph{not} a multiple of
$T_{ab}$ at some $m'\in V$ then (\ref{50}) shows, since $T(m')\neq
0$, that $\rho(m')=\sigma(m')=0$ and then (\ref{43}) gives the
contradiction that $\lambda(m')=0$. So $p_{a;b}$ is a
\emph{non-zero} multiple of $T_{ab}$ at each point of $V$. It
follows that the covector f\/ields associated with both $l$ and $p$
are closed 1-forms on $V$ and hence are gradients on (a possibly
reduced) $V$, say of the (smooth) functions $u$ and $v$,
respectively. In addition, $l$ and $p$ are easily checked to have
zero Lie bracket and so span an integrable distribution (again on a
possibly reduced $V$). It can then be shown by using standard
techniques that if one chooses $x$ and $y$ coordinates on that part
of the 2-dimensional submanifold $u=v=0$ contained in (a possibly
reduced) $V$ and chooses the other coordinates as parameters along
the integral curves of~$l$ and~$p$ (with~$x$ and~$y$ constant along
these curves) the parameter coordinates may be taken as $u$ and $v$
and $l=\partial /\partial v$ and $p=\partial /\partial u$. Then with
the coordinates  $u$, $v$, $x$, $y$ the metric $g$ takes the form
\begin{equation}\label{51}
ds^{2}=2dudv+g_{\alpha \beta}dx^{\alpha}dx^{\beta} \qquad (g_{\alpha \beta}=T_{\alpha \beta}), \end{equation}
where Greek indices take the values 3 and 4 and the second (bracketed) equation in (\ref{47}) has been used.

Now suppose that the constant $c\neq 0$. Then (\ref{46}) shows that
one may use the translational freedom in the coordinate $u$ to
arrange that $\rho=cu$. If $c=0$ then the f\/irst equation in~(\ref{47}) together with the non-vanishing of $p_{a;b}$ on $V$
implies $\rho=0$ on $V$ and the second equation in~(\ref{46}) shows
that $\sigma$ is a non-zero constant on~$V$. Then~(\ref{49}) and~(\ref{50}) reveal that $\delta$ is a linear function of $u$ which is
nowhere zero on~$V$. A translation of the $u$ coordinate then
ensures that $\delta=\sigma u$ (and $u$ is nowhere zero) on $V$.
Thus whether $c$ is zero or not, one can achieve the result that
$up_{a;b}=T_{ab}$. Denoting the Christof\/fel symbols arising from
$\nabla$ by $\Gamma$ and noting that since $l=\partial /\partial v$
is a Killing vector f\/ield on~$V$, the functions $g_{\alpha \beta}
(=T_{\alpha \beta}$) are independent of $v$, one calculates that
$-u\Gamma^{2}_{\alpha \beta}=g_{\alpha \beta}$. This, together with
$up_{a;b}=T_{ab}$ and an integration then leads to~$g_{\alpha
\beta}=u^{2}h_{\alpha \beta}$ for functions $h_{\alpha \beta}$ which
are independent of $u$ and $v$. Thus, on $V$, the original metric is
\begin{equation}\label{52}
ds^{2}=2dudv+u^{2}h_{\alpha \beta}\big(x^{3},x^{4}\big)dx^{\alpha}dx^{\beta}.
\end{equation}

The next step is, recalling the coordinate scalings above, to write
$\rho=cu$ and, from (\ref{46}), $\sigma=cv+e_{1}$ and substitute
into (\ref{49}) to get $\beta=cv^{2}+2e_{1}v+e_{2}$,
$\gamma=cu^{2}+e_{3}$ and $\delta=e_{1}u+cuv+e_{4}$ for constants
$e_{1}$, $e_{2}$, $e_{3}$ and $e_{4}$. Then (\ref{47}) and
(\ref{50}) act as consistency checks on these constants and give
$e_{3}=e_{4}=0$. The equation (\ref{49}) now gives the Sinyukov
tensor in (\ref{45}) as \begin{equation}\label{53}
a_{ab}=\phi\left(g_{ab}+(cv^{2}
+2e_{1}v+e_{2})l_{a}l_{b}+cu^{2}p_{a}p_{b}+(e_{1}u+cuv)(l_{a}p_{b}+p_{a}l_{b})\right).
\end{equation} This equation can then be inverted following the
techniques described in Section~\ref{section5}. First, with the above values for
$\rho$ and $\sigma$ one f\/inds $\lambda$ from (\ref{43}) and, from
(\ref{53}), $\psi$ can be calculated from the equation
$\psi_{a}=-a^{-1}_{ab}\lambda^{b}$ given in Section~\ref{section5}. Then the
potential, $\chi$, of the 1-form $\psi$ can be found. Finally the
required metric is then given, on $V$, by $g'_{ab}=e^{2\chi}a_{ab}$.
This calculation was done using Maple and results in
\begin{equation}\label{54}
\psi=d\chi, \qquad \chi=\tfrac{1}{2}\ln F,\qquad F=\kappa^4\left[1+2cuv+2e_{1}u+(e_{1}^{2}-ce_{2})u^{2}\right]^{-1} \end{equation}
and for the metric $g'$, projectively related to $g$,
\begin{gather} ds'^{2}=\kappa Fg-\kappa^{-3}F^{2}\big[(cv^{2}+2e_{1}v+e_{2})du^{2}+cu^{2}dv^{2} \nonumber\\
\phantom{ds'^{2}=}{}
+2u(cv+e_{1}+(e_{1}^{2}-ce_{2})u)dudv\big] \label{55}
\end{gather} for a positive constant $\kappa=\phi^{-1}$.

It is noted here that if $\lambda_{a}(m)$ is proportional to $l(m)$
then it follows from (\ref{43}) that $\rho(m)=0$ and from~(\ref{47})
that $c=0$. Then (\ref{46}) reveals that $\rho=0$ and that $\sigma$
is constant on $V$ and so $\lambda_{a}$ is a constant multiple of
$l_{a}=g_{ab}l^b$ on~$V$. Conversely, since $\nabla \lambda=cg$, this last
condition implies that $c=0$. Thus the special case $c=0$
corresponds to the condition that $\lambda_{a}$ is a~(non-zero)
constant multiple of $l_{a}$ on $V$ and then~(\ref{53}), (\ref{54})
and (\ref{55}) simplify. For this special case it follows from
(\ref{54}) that $\psi_{a}=r(u)l_{a}$ on $V$ for some (smooth)
function $r$ on $V$. Thus, for any metric $g'$, projectively related
to $g$ on $V$ and with Levi-Civita connection $\nabla'$, one f\/inds
from (\ref{E12}) that $l_{a|b}=-2r(u)l_{a}l_{b}$, where a stroke
denotes a $\nabla'$-covariant derivative. Then the nowhere-zero vector f\/ield $l'\equiv
e^{\int 2r}l$ on $V$ is null with respect to $g'$ (from~(\ref{55}) with $c=0$)
and covariantly constant with respect to~$\nabla'$. It is, in fact,
true that $(V,g')$ (that is, $\nabla'$) is also of holonomy type~$R_{11}$.
This follows from Table~\ref{table1} since $(V,g')$ admits a nowhere
zero null covariantly constant vector f\/ield and so must be of
holonomy type $R_{3}$, $R_{4}$, $R_{8}$ or $R_{11}$ and the f\/irst
three of these are ruled out because they would, from Theorems~\ref{theorem2} and~\ref{theorem5}, force $\nabla=\nabla'$. [That $\nabla'$ cannot be f\/lat over some
non-empty open subset $V'$ of $V$ follows since $g'$ would then be
of constant curvature and hence so would $g$~\cite{15},
contradicting the curvature class $\mathbf{D}$ condition on
$(V,g)$.] The space-time $(V,g')$ is also of curvature class
$\mathbf{D}$ on $V$ since, otherwise, it would be of curvature class
$\mathbf{C}$ over some non-empty open subset~$V'$ of~$V$. In this
case, Theorem~\ref{theorem4} would show that $\lambda=0$ on~$V'$ and a
contradiction follows.

Now let $(M,g)$ be a space-time of holonomy type $R_{10}$ or
$R_{13}$ which is (non-f\/lat and) of curvature class $\mathbf{D}$,
that is the curvature rank equals 1 over some open dense subset $U$
of $M$ and $\rm Riem$ vanishes on $M\setminus U$. Let $g'$ be a metric
on some connected open neighbourhood $V\subset U$ of $m\in U$ which
is projectively related to~$g$ on~$V$. On $V$, $\rm Riem$ takes the form
(\ref{E5}) for some nowhere zero bivector $F$ on $V$. If the
holonomy type is $R_{13}$ this bivector $F$ is necessarily spacelike
and so the subset $\mathbf{D_{n}}$ in the above decomposition of $M$ is
empty. Then $V$ may be chosen to admit a~covariantly constant, unit
timelike vector f\/ield $u$ (so that, from the Ricci identity,
$R^{a}{}_{bcd}u^{d}=0$ on~$V$) and a unit (smooth) spacelike vector
f\/ield $w$ which is orthogonal to $u$ on~$V$ and also satisf\/ies
$R^{a}{}_{bcd}w^{d}=0$ on~$V$. On $V$ one has the useful relations
$u^{a}u_{a;b}=w^{a}w_{a;b}=u^{a}w_{a;b}=0$. If, however, the
holonomy type of $(M,g)$ is $R_{10}$, $V$ may be chosen so as to
admit a covariantly constant unit spacelike vector f\/ield (which for
later convenience will also be labelled $u$) and which satisf\/ies
$R^{a}{}_{bcd}u^{d}=0$ and another nowhere-zero vector f\/ield
(similarly labelled $w$) which is orthogonal to $u$ and satisf\/ies
$R^{a}{}_{bcd}w^{d}=0$ on $V$ but whose nature (spacelike, timelike
or null) may vary. It follows that in the $R_{10}$ case, each of
$\mathbf{D_{n}}$ and $\mathbf{D_{nn}}$ may be non-empty. These remarks permit the
totality of possibilities for the types $R_{10}$ and $R_{13}$, with
one exception, to be taken together (the exceptional case being for
type~$R_{10}$ and the subset $\mathbf{D_{n}}$).

The proof is, in fact, very similar to that already given for the
holonomy type $R_{11}$ and so it will be described only brief\/ly. For
$m\in U$ choose an open, connected and simply connected
neighbourhood $V$ of $m$. Then one has, from Lemma~\ref{lemma1},
$\lambda_{a;b}=cg_{ab}$ for some constant $c$ and so, as before,
\begin{equation}\label{56}
\lambda_{a}=\sigma u_{a}+\rho w_{a}, \qquad \rho w_{a;b}+w_{a}\rho _{,b}+u_{a}\sigma _{,b}=cg_{ab} \end{equation} and
\begin{equation}\label{57}
a_{ab}=\phi g_{ab}+\beta u_{a}u_{b}+\gamma w_{a}w_{b}+\delta (u_{a}w_{b}+w_{a}u_{b})
\end{equation} for functions $\sigma$, $\rho$, $\phi$, $\beta$, $\gamma$ and $\delta$ on $V$.
Whichever subcase is chosen an argument similar to that given in the
holonomy $R_{11}$ case leads to the conclusion that $\phi$ is
constant.

Now suppose the holonomy type is $R_{10}$, that $\mathrm{int}\,(M_{2}\cap
\mathbf{D_{n}})$ is non-empty and with $V$ as usual choose $m\in V\subset
\mathrm{int}\,(M_{2}\cap \mathbf{D_{n}})$. On $V$, $u$ is unit spacelike and covariantly
constant and $w$ is null. Then the second of (\ref{56}), on
contraction with $w^{a}$ and $u^{a}$, gives $c=0$, $\sigma$ is
constant on $V$ and $(\rho w)_{a;b}=0$. Thus $\lambda$ is
covariantly constant on $V$ and $\rho(m)=0$. [If $\rho(m)\neq 0$
then $\rho$ would not vanish over some neighbourhood of $m$ and on
that neighbourhood $\rho w$ would be a nowhere-zero, nowhere-null
covariantly constant vector f\/ield on $V$ with $\rho w(m)$ and $u(m)$
independent members of $T_{m}M$ at each $m\in V$ . This would reduce
the holonomy type in this neighbourhood to $R_{3}$ and then, from
Theorem~\ref{theorem2}, one f\/inds that $\lambda$ vanishes on this neighbourhood
of $m$ contradicting the fact that $m\in M_{2}$.] Since $\lambda$
(and $u$) are covariantly constant on $V$ with $\lambda$
proportional to $u$ at $m$ it follows that $\lambda$ is proportional
to $u$ on $V$ and hence, from~(\ref{56}), that~$\rho$ is zero on
$V$. Finally a substitution of~(\ref{57}) into~(\ref{E18}) and a
contraction with $u^{a}w^{b}$ gives $\sigma=0$ on $V$ and hence,
from~(\ref{56}) that $\lambda\equiv 0$ on $V$. This contradiction to
the statement $m\in M_{2}$ shows that $\mathrm{int}(M_{2}\cap
\mathbf{D_{n}})=\varnothing$.

The other cases can be handled together. One assumes that $M_{2}\cap
\mathbf{D_{nn}}$ is non-empty and chooses $m\in (M_{2}\cap \mathbf{D_{nn}})$ and the
usual neighbourhood $V$ such that $m\in V\subset \mathbf{D_{nn}}$. To deal
with both cases simultaneously one writes, for the vector f\/ields $u$
and $w$ on $V$, $u^{a}u_{a}=\epsilon _{1}$, $w^{a}w_{a}=\epsilon
_{2}$ (and $u^{a}w_{a}=0$) with $\epsilon _{1}$ and $\epsilon _{2}$
equal to either $\pm 1$. Thus the $R_{13}$ type requires $\epsilon
_{1}=-1$ and $\epsilon _{2}=1$ whilst the $R_{10}$ type requires
$\epsilon _{1}=1$ and $\epsilon _{2}=\pm 1$. Equations (\ref{56})
and (\ref{57}) still hold and contractions of the second in
(\ref{56}) with $u^{a}$ and $w^{a}$ give
\begin{equation}\label{58}
\epsilon _{1}\sigma _{,a}=cu_{a}, \qquad \epsilon _{2}\rho _{,a}=cw_{a}, \qquad \rho w_{a;b}=cT_{ab}\qquad (T_{ab}=g_{ab}- \epsilon _{1}u_{a}u_{b}-\epsilon _{2}w_{a}w_{b}). \end{equation}
Then substituting (\ref{57}) into (\ref{E18}) and contracting successively with $u^{a}u^{b}$, $w^{a}w^{b}$ and $u^{a}w^{b}$ and then with $u^{a}$ and $w^{a}$ gives
\begin{equation}\label{59}
2\epsilon _{1}\sigma u_{a}=\beta _{,a},\qquad 2\epsilon _{2}\rho w_{a}=\gamma _{,a}, \qquad \epsilon _{2} \rho u_{a}+\epsilon _{1} \sigma w_{a}=\epsilon _{1}\epsilon _{2} \delta _{,a}
\end{equation}
and
\begin{equation}\label{60}
\gamma w_{a;b}=\rho T_{ab},\qquad \delta w_{a;b}=\sigma T_{ab}.
\end{equation}
As before one can now show, since $\lambda\in M_{2}$, that $\lambda$
and $w_{a;b}$ are nowhere zero on $V$ and hence $V\subset M_{2}\cap
\mathbf{D_{nn}}$ (so that $M_{2}\cap \mathbf{D_{nn}}$ is open in $M$). Also, on $V$,
$w_{a;b}$ is a non-zero multiple of~$T_{ab}$. Thus, $V$ may be
chosen so that the 1-forms associated with $u$ and $w$ through $g$
are global gradients on $V$, say $u=\epsilon _{1}dt$ and $w=\epsilon
_{2}dz$ for functions $t$ and $z$ on $V$. So reducing $V$, if
necessary, to a coordinate domain with coordinates $t$, $z$, $x^{3}$,
$x^{4}$ such that $g$ takes the form
\begin{equation}\label{61}
ds^{2}=\epsilon_{1}dt^{2}+\epsilon_{2}dz^{2}+z^{2}h_{\alpha \beta}(x^{3}, x^{4})dx^{\alpha} dx^{\beta},
\end{equation}
where the functions $h_{\alpha \beta}$ are independent of $t$ and $z$. The equations (\ref{56}) and (\ref{59}) may be integrated to give
\begin{equation}\label{D1}
a_{ab}=\phi\left(g_{ab}+(ct^2+2c_2t+c_3)u_au_b+cz^2w_aw_b+(ctz+c_2z)(u_aw_b+w_au_b)\right),
\end{equation}
where $c_2$ and $c_3$ are constants. Inverting the pair
$(a,\lambda)$ to obtain the corresponding pair $(g',\psi)$ leads to
\begin{gather}
\psi=d\chi, \qquad \chi=\tfrac{1}{2}\ln F, \nonumber\\
\label{D2} F=\kappa^4\left[1+\epsilon_2(c+\epsilon_1(c_3c-c_2^2))z^2
+\epsilon_1(ct^2+2c_2t+c_3)\right]^{-1}
\end{gather}
corresponding to the metric $g'$ projectively related to $g$ and
given by
\begin{equation}\label{D3}
g'=\kappa Fg-\kappa^{-3}F^2\left(\begin{array}{l}(ct^2+2c_2t+c_3+\epsilon_2(cc_3-c_2^2)z^2
)dt^2\\
+(c+\epsilon_1(cc_3-c_2^2))z^2dz^2
+2\epsilon_1\epsilon_2(ct+c_2)zdtdz\end{array}\right),
\end{equation}
where $\kappa=\phi^{-1}$ is a positive constant.

In summary, for space-times $(M,g)$ of holonomy types $R_{10}$,
$R_{11}$ and $R_{13}$ and which are of curvature class $\mathbf{D}$
over some open dense subset $U$ of $M$, $\mathrm{int} (M_{2} \cap
\mathbf{D_{n}})=\varnothing$ and one may identify an open dense subset of $M$
such that each point $m$ of this subset lies in a neighbourhood
where either $\lambda=0$ (when $m\in M_{1}$) or where $\lambda$ is
nowhere zero and in which the situation can be completely resolved
(in the open subset $M_{2}\cap \mathbf{D_{nn}}$). It is interesting to remark
at this point that more can be said in the cases when $\rm Riem$ is
\emph{nowhere zero} on $M$ and the \emph{nature} (timelike,
spacelike or null) of the bivector $F$ in (\ref{E5}) is constant on
$M$. To see this, suppose that such is the case and suppose that $F$
is non-null (and hence always timelike or spacelike on~$M$). Then, in the previous notation, $M=M_{1}\cup
M_{2}=\mathbf{D_{nn}}$ with $M_{1}$ open in $M$. But then the previous work
shows that $M_{2}=M_{2}\cap \mathbf{D_{nn}}$ is also open in $M$. Thus, since
$M$ is connected, one of $M_{1}$ and $M_{2}$ is empty and so either
$M=M_{1}$ or $M=M_{2}$. It follows that, in this case, either
$\lambda=0$ on $M$ or $\lambda$ is nowhere vanishing on $M$. In the
similar situation, but with $F$ null, $M=\mathbf{D_{n}}$ and the work above
shows that $\lambda=0$ on $M$.

Although the class $\mathbf{D}$ assumption for these holonomy types
was made over an open dense subset of $M$ it is now seen that this
is, in fact, a complete solution for holonomy types $R_{10}$,
$R_{11}$ and $R_{13}$. For if $(M,g)$ is of any of these holonomy
types, then either it is of curvature class $\mathbf{C}$ over an
open dense subset of $M$, in which case, $\nabla=\nabla'$ on $M$, or
the curvature class is $\mathbf{D}$ over some subset $E$ of $M$. Now
the subset $J$ of $M$ over which $\rm Riem$ vanishes is closed and
satisf\/ies $\mathrm{int}\,J=\varnothing$. Also, the subset $E\cup J$ is closed
(by the rank theorem on $\rm Riem$ since $M\setminus (E\cup J)$ is
exactly the open subset of $M$ on which $\rm Riem$ has curvature rank $>
1$) and if $\mathrm{int}E=\varnothing$ it is easily checked that $\mathrm{int}\, (E\cup
J)=\varnothing$ and so the curvature class is $\mathbf{C}$ over the
open dense subset $M\setminus (E\cup J)$ of $M$. Then Theorem~\ref{theorem4}
gives $\nabla=\nabla'$ on $M$. If $\mathrm{int}E\neq \varnothing$ let $E'=\mathrm{int}\,
E$ and proceed as above on the (connected) components of $E'$ with
metric $g$.

Where a space-time $(M,g)$ has holonomy type $R_{9}$ or $R_{14}$ the
situation turns out to be even more complicated. These cases will be
discussed in detail elsewhere~\cite{26}, but, for
completeness, non-trivial (that is, $\lambda$ not identically zero)
examples of projective relatedness in such space-times will now be
discussed brief\/ly.

Consider a spacetime $(M,g)$ and a point $m\in M$ for which there
exists a coordinate neighbourhood $U$ of $m$ with coordinates
$u$, $v$, $x$, $y$ ($v> 0$) such that $g$ is given on $U$ by
\begin{equation}\label{A.1}
ds^2=2dudv+b(u)\sqrt{v}du^2+u^2e^{f(x,y)}\big(dx^2+dy^2\big),
\end{equation}
where $b(u)$ and $f(x,y)$ are arbitrary smooth functions on $U$. Now
the 1-form f\/ield $l=du$  may be shown to be recurrent on $U$, and is
covariantly constant on $U$ if and only if $b(u)\equiv 0$ on $U$. On
any connected open subset of $U$  where $b(u)=0$,~(\ref{A.1}) is
locally of the form~(\ref{52}) and of curvature class~$\mathbf{D}$.
However, for~(\ref{A.1}) to be of holonomy type $R_{9}$ or $R_{14}$,
$b(u)\neq0$ at least at one point $m\in M$, and hence over some open
neighbourhood $V$ of $m$. Assuming that $b(u)$ is nowhere zero on
$V$ the curvature class is $\mathbf{A}$ on $V$ and after somewhat lengthy
calculations it may be shown that~(\ref{E18}) may be solved
non-trivially to give
\begin{equation}\label{A.2}
a=\phi\big(g+2\xi ududv+\xi\big(2v+ub(u)\sqrt{v}\big)du^2\big),
\end{equation}
where $\phi,\xi\in\mathbb{R}$. Inverting the pair $(a,\lambda)$ to
obtain the corresponding pair $(g',\psi)$ leads to
\begin{equation}\label{A.3}
\psi=d\chi, \qquad \chi=\tfrac{1}{2}\ln F, \qquad F=\kappa^4\left[1+\xi u\right]^{-2}
\end{equation}
with $\kappa=\phi^{-1}$, corresponding to a metric $g'$,
projectively related to (\ref{A.1}), and given on $V$ by
\begin{equation}\label{A.4}
g'=\kappa F g-\kappa^{-3}F^2\left[2\xi u(1+\xi u)dudv+\xi(u(1+\xi u)b(u)\sqrt{v}+2v)du^2\right].
\end{equation}
If $f(x,y)$ satisf\/ies $\tfrac{\partial^2f}{\partial
x^2}+\tfrac{\partial^2f}{\partial y^2}=0$ everywhere on $V$ then an
application of the Ambrose--Singer theorem \cite{13} shows
(\ref{A.1}) has holonomy type~$R_9$, and otherwise has holonomy type
$R_{14}$ on $V$.

It is remarked that examples are known of non-trivially projectively
related space-times of the Friedmann--Robertson--Walker--Lemaitre
(FRWL) type and which are of holonomy type~$R_{15}$~\cite{27}. Other examples of $R_{15}$ metrics with similar properties
will be discussed in~\cite{26}.

\section{Conclusions}\label{section8}

This paper studied the question of, given a space-time $(M,g)$, what can be said about those other metrics $g'$ on $M$ (or some open submanifold of $M$) which are
projectively related to~$g$? The question was inspired from both a
geometrical point of view and from a physical standpoint through the
principle of equivalence in Einstein's general theory of relativity.
If $(M,g)$ is an Einstein space (and this includes the important
vacuum metrics of general relativity) the problem is already
resolved, as explained in Section~\ref{section5}. This paper then proceeded to
extend this by studying the situation through the various possible
holonomy groups that are possible for a~Lorentz manifold. It was
shown that, by using certain algebraic properties of the curvature
tensor, much of this problem is tractable and an essentially
complete solution was obtained for all the holonomy types except
$R_{9}$, $R_{14}$ and $R_{15}$. For these latter types, examples of
non-trivially projectively related space-times were shown to exist
(and further details will be given elsewhere). Amongst all of these
are several well-known, physically important solutions of Einstein's
f\/ield equations.

\pdfbookmark[1]{References}{ref}
\LastPageEnding


\begin{thebibliography}{99}

\footnotesize\itemsep=0pt

\bibitem{1}
Geroch R.P.,
Spinor structure of space-times in general relativity,
{\it J. Math. Phys.} {\bf 9} (1968), 1739--1744.

\bibitem{2}
Hall G.S.,
Symmetries and curvature structure in general relativity, {\it World Scientific Lecture Notes in Physics}, Vol.~46, World Scientif\/ic Publishing Co., Inc., River Edge, NJ, 2004.

\bibitem{3}
Hall G.S., McIntosh C.G.B.,
Algebraic determination of the metric from the curvature in general relativity,
{\it Internat. J. Theoret. Phys.}  {\bf 22} (1983),   469--476.

\bibitem{4}
Hall G.S.,
Curvature collineations and the determination of the metric from the curvature in general relativity,
{\it Gen. Relativity Gravitation}  {\bf 15}  (1983),  581--589.

\bibitem{5}
McIntosh C.G.B., Halford W.D.,
Determination of the metric tensor from components of the Riemann tensor,
{\it J. Phys. A: Math. Gen.}  {\bf 14} (1981),  2331--2338.

\bibitem{6}
Hall G.S., Lonie D.P.,
On the compatibility of Lorentz metrics with linear connections on four-dimensional manifolds,
{\it J. Phys. A: Math. Gen.}  {\bf 39}  (2006), 2995--3010,
\href{http://arxiv.org/abs/gr-qc/0509067}{gr-qc/0509067}.

\bibitem{7}
Kobayashi S., Nomizu K.,
Foundations of dif\/ferential geometry, Vol.~1, Interscience Publishers, New York,  1963.

\bibitem{8}
Schell J.F.,
Classif\/ication of four-dimensional Riemannian spaces,
{\it J. Math. Phys.} {\bf 2} (1961),  202--206.

\bibitem{25}
Hall G.S., Lonie D.P.,
Holonomy groups and spacetimes,
{\it Classical Quantum Gravity} {\bf 17} (2000), 1369--1382,
\href{http://arxiv.org/abs/gr-qc/0310076}{gr-qc/0310076}.

\bibitem{9}
Hall G.S.,
Covariantly constant tensors and holonomy structure in general relativity,
{\it J. Math. Phys.} {\bf 32} (1991), 181--187.

\bibitem{10}
Besse A.,
Einstein manifolds,  Springer-Verlag, Berlin, 1987.

\bibitem{11}
Wu H.,
On the de Rham decomposition theorem,
{\it Illinois J. Math.} {\bf 8} (1964),  291--311.

\bibitem{12}
Hall G.S., Kay W.,
Holonomy groups in general relativity,
{\it J. Math. Phys.}  {\bf 29} (1988), 428--432.

\bibitem{13}
Ambrose W., Singer I.M.,
A theorem on holonomy,
{\it Trans. Amer. Math. Soc.} {\bf 75} (1953), 428--443.

\bibitem{14}
Hall G.S.,
Connections and symmetries in space-times,
{\it Gen. Relativity Gravitation}  {\bf 20}  (1988), 399--406.

\bibitem{15}
Eisenhart L.P.,
Riemannian geometry,  Princeton University Press, Princeton, 1966.

\bibitem{16}
Thomas T.Y.,
Dif\/ferential invariants of generalised spaces,  Cambridge,  1934.

\bibitem{17}
Weyl H.,
Zur Inf\/initesimalgeometrie: Einordnung der projectiven und der konformen Auf\/fassung,
{\it G\"ott. Nachr.} (1921), 99--112.

\bibitem{18}
Petrov A.Z.,
Einstein spaces,  Pergamon Press, Oxford~-- Edinburgh~--New York, 1969.

\bibitem{19}
Sinyukov N.S.,
Geodesic mappings of Riemannian spaces, Nauka, Moscow, 1979 (in Russian).

\bibitem{20}
Mikes J., Hinterleitner I., Kiosak V.A.,
On the theory of geodesic mappings of Einstein spaces and their gene\-ralizations,  in The Albert Einstein Centenary International Conference, Editors J.-M.~Alini and A.~Fuzfa, {\it AIP Conf. Proc.}, Vol. 861, American Institute of Physics, 2006, 428--435.

\bibitem{21}
Hall G.S., Lonie D.P.,
The principle of equivalence and projective structure in spacetimes,
{\it Classical Quantum Gravity}  {\bf 24} (2007), 3617--3636,
\href{http://arxiv.org/abs/gr-qc/0703104}{gr-qc/0703104}.

\bibitem{22}
Kiosak V., Matveev V.A.,
Complete Einstein metrics are geodesically rigid,
{\it Comm. Math. Phys.} {\bf 289} (2009), 383--400,
\href{http://arxiv.org/abs/0806.3169}{arXiv:0806.3169}.

\bibitem{23}
Hall G.S., Lonie D.P.,
Projective equivalence of Einstein spaces in general relativity,
{\it Classical Quantum Gravity} {\bf 26} (2009), 125009, 10~pages.

\bibitem{26}
Hall G.S., Lonie D.P.,
Holonomy and projective structure in space-times, Preprint, University of Aberdeen, 2009.

\bibitem{27}
Hall G.S., Lonie D.P.,
The principle of equivalence and cosmological metrics,
{\it J. Math. Phys.} {\bf 49} (2008), 022502, 13~pages.

\end{thebibliography}
\end{document}